\documentclass[12pt,reqno]{amsart}
\usepackage[T1]{fontenc}
\usepackage[nomath]{lmodern}
\usepackage[margin=1in]{geometry}
\usepackage{amsmath,amssymb,amsthm}
\usepackage[hidelinks]{hyperref}

\numberwithin{equation}{section}
\allowdisplaybreaks[1]
\newtheorem{theorem}{Theorem}[section]
\newtheorem{proposition}[theorem]{Proposition}
\newtheorem{lemma}[theorem]{Lemma}
\newtheorem{corollary}[theorem]{Corollary}
\theoremstyle{remark}
\newtheorem{remark}[theorem]{Remark}

\newcommand{\R}{\mathbb R}
\newcommand{\C}{\mathbb C}
\newcommand{\Sph}{\mathbb S}
\newcommand{\dbar}{\bar\partial}
\newcommand{\norm}[1]{\lVert#1\rVert}
\newcommand{\ip}[2]{\langle#1,#2\rangle}
\newcommand{\Div}{\operatorname{div}}
\newcommand{\Divmu}{\operatorname{div}_\mu}
\newcommand{\T}{\mathcal T}
\newcommand{\Or}{\mathcal O}
\newcommand{\X}{\mathcal X}
\newcommand{\eps}{\varepsilon}
\renewcommand{\Re}{\operatorname{Re}}
\renewcommand{\Im}{\operatorname{Im}}

\title[Sharp Dirac--Sobolev inequalities in low dimensions]{Sharp Dirac--Sobolev inequalities in low dimensions}
\author[J. Seo]{Juneyoung Seo}
\address{Innovation Center for Atomic Science, DAU G-LAMP Project Group,
Dong-A University, Busan 49315, Republic of Korea}
\email{juneys.dau@gmail.com}
\thanks{This research was supported by Global - Learning \& Academic research institution for Master's$\cdot$PhD students, and Postdocs (LAMP) Program of the National Research Foundation of Korea (NRF) grant funded by the Ministry of Education (RS-2025-25440216).}

\begin{document}

\begin{abstract}
In dimensions $d=2,3,4$, the sharp constant in the Dirac--Sobolev inequality on $\R^d$
\[
 \norm{D_d\psi}_{p_d}\ge S_d\norm{\psi}_{q_d},
 \qquad p_d=\frac{2d}{d+1},\quad q_d=\frac{2d}{d-1},
\]
is the best Sobolev constant for $\dot W^{1,p_d}(\R^d)\hookrightarrow L^{q_d}(\R^d)$. We classify all optimizers. In the plane, the area formula for an anti-holomorphic rational map associated with an optimizer determines the constant. In dimensions three and four, we use a conformal metric to estimate volume and a second integral identity to control the $L^{q_d}$ norm of solutions to the coupled Euler--Lagrange system. As consequences we obtain the sharp planar Cauchy transform norm $2/\pi$ and the sharp gradient estimate $\norm{\nabla u}_{q_d}\le S_d^{-1}\norm{\Delta u}_{p_d}$.
\end{abstract}

\maketitle

\section{Introduction}

\subsection{Problem and main result}
Let $d\ge2$ and $N_d=2^{\lfloor d/2\rfloor}$. Fix skew-Hermitian matrices $c_1,\dots,c_d\in\C^{N_d\times N_d}$ satisfying
\[
 c_j^*=-c_j,\qquad c_jc_k+c_kc_j=-2\delta_{jk}I_{N_d},
\]
and acting irreducibly on $\C^{N_d}$. For $x\in\R^d$ and $z,w\in\C^{N_d}$, we write $c(x)=\sum_{j=1}^d x_jc_j$, $|z|=\bigl(\sum_{\alpha=1}^{N_d}|z_\alpha|^2\bigr)^{1/2}$, and $\ip{z}{w}=\Re\sum_{\alpha=1}^{N_d}z_\alpha\overline{w_\alpha}$. The Dirac operator $D_d=\sum_{j=1}^d c_j\partial_j$ on $\R^d$ is formally self-adjoint with $D_d^2=-\Delta$. Set
\[
 p_d=\frac{2d}{d+1},\qquad q_d=\frac{2d}{d-1},
\]
so that $1/p_d+1/q_d=1$ and $1/p_d-1/q_d=1/d$. On the space
\[
 \X_d=\bigl\{\psi\in L^{q_d}(\R^d;\C^{N_d}):D_d\psi\in L^{p_d}(\R^d;\C^{N_d})\text{ distributionally}\bigr\},
\]
we consider the minimization problem
\[
 \inf_{\psi\in\X_d\setminus\{0\}}Q_d(\psi),\qquad Q_d(\psi)=\frac{\norm{D_d\psi}_{p_d}}{\norm{\psi}_{q_d}}.
\]
The infimum is attained for every $d\ge2$ by concentration compactness (cf.\ \cite{FL-exist} for $d=3$).

If a real magnetic potential $A\in L^d(\R^d;\R^d)$ admits a zero mode $0\ne\psi\in L^{q_d}$ of $D_d-ic(A)$, then H\"older's inequality gives $Q_d(\psi)\le\|A\|_{L^d}$. Thus a lower bound on $Q_d$ gives a necessary lower bound on the norm of such a potential. The zero-mode problem arises in the study of the stability of matter in magnetic fields \cite{FL22,FL,LLS} (see also \cite{FL-exist,Mederski,WZzero,WZcurl}).

If $\psi=f\eta$ for a real-valued function $f$ and a fixed unit spinor $\eta\in\C^{N_d}$, then $|f\eta|=|f|$ and $|D_d(f\eta)|=|\nabla f|$. Thus $\inf Q_d$ is bounded above by the best Sobolev constant \cite{Aubin,Talenti}
\begin{equation}\label{eq:Sd}
 S_d=\inf_{0\ne f\in\dot W^{1,p_d}(\R^d)}\frac{\norm{\nabla f}_{p_d}}{\norm{f}_{q_d}}
 =d\biggl(\frac{|\Sph^{d-1}|}{q_d}\,\mathrm B\Bigl(\frac{d-1}{2},\frac{d+1}{2}\Bigr)\biggr)^{1/d},
\end{equation}
where $\dot W^{1,p_d}(\R^d)=\{f\in L^{q_d}(\R^d):\nabla f\in L^{p_d}(\R^d)\}$ and $|\Sph^{d-1}|$ denotes the area of the unit sphere in $\R^d$. In $\dot W^{1,p_d}(\R^d)$, the constant $S_d$ is attained by the Aubin--Talenti profile
\[
 f_d(r)=\bigl(1+r^{q_d}\bigr)^{-(d-1)/2}\qquad(r\ge0).
\]

Dirac--Sobolev estimates in dimension three were studied in \cite{BES,BEU}. Frank and Loss posed the sharp problem in that dimension \cite{FL22}. In every dimension $d\ge2$, stereographic projections of Killing spinors on $\Sph^d$ are critical points with quotient $(d/2)|\Sph^d|^{1/d}$ \cite{FL} (see also \cite{BMW,Hijazi,Xu}). Frank and Loss conjectured that this value was sharp \cite{FL-exist,FL}. Wang and Zhang~\cite[Theorem~5.6]{WZ} disproved the conjecture and showed that the corresponding quotient on $\Sph^d$ has Morse index $d+1$ at each element of the conformal Killing-spinor orbit with respect to smooth variations. Bieganowski and Mederski showed that the constant-direction trial spinor $f_d(|x|)\eta$ has quotient $S_d<(d/2)|\Sph^d|^{1/d}$ \cite{BM}.

When the spinor direction varies with $x$, the components of $D_d\psi$ couple and $|D_d\psi|$ can be pointwise smaller than $|\nabla\psi|$. On $\C\cong\R^2$, with $2\partial=\partial_1-i\partial_2$ and $2\dbar=\partial_1+i\partial_2$, any smooth $u\colon\C\to\C$ satisfies
\[
 |2\dbar u|^2=|\nabla u|^2-2\det\nabla(\Re u,\Im u).
\]
The integrated identity $\norm{D_d\psi}_2=\norm{\nabla\psi}_2$ does not give a comparison of these norms at the exponent $p_d<2$.

For dimensions $d=2,3,4$, we prove $\inf Q_d=S_d$. All optimizers are obtained from the constant-direction profiles by translation, dilation, and the spinorial Kelvin transform.

\begin{theorem}\label{thm}
Let $d\in\{2,3,4\}$. For every $\psi\in\X_d$,
\begin{equation}\label{eq:sharp}
 \norm{D_d\psi}_{p_d}\ge S_d\,\norm{\psi}_{q_d},
\end{equation}
and the constant $S_d$ is sharp, with $S_2=\pi$, $S_3=(18\pi)^{1/3}$, and $S_4=(12\pi^3)^{1/4}$. Equality holds in \eqref{eq:sharp} for $\psi\ne0$ if and only if, for some $\eta\in\C^{N_d}\setminus\{0\}$,
\begin{equation}\label{eq:optimizers}
 \psi(x)=\frac{\eta}{\bigl(\sigma^{q_d}+|x-x_0|^{q_d}\bigr)^{(d-1)/2}}
 \qquad\text{or}\qquad
 \psi(x)=\frac{c(x-x_1)\eta}{\bigl(|x-x_1|^{q_d}+|\alpha x+\beta|^{q_d}\bigr)^{(d-1)/2}},
\end{equation}
where $\sigma>0$ and $x_0\in\R^d$ in the first family, and $x_1\in\R^d$, $\alpha\ge0$, and $\beta\in\R^d$ with $\alpha x_1+\beta\ne0$ in the second.
\end{theorem}

For $d\ge5$, the sign estimates in this paper do not determine the sharp constant.

The spinorial Kelvin transform \cite[Section~4]{BF}
\begin{equation}\label{eq:kelvin-spinor}
 \psi_{\mathcal K}(x)=|x|^{-d}c(x)\,\psi\bigl(x/|x|^2\bigr),\qquad (\psi_{\mathcal K})_{\mathcal K}=-\psi,
\end{equation}
preserves $\norm{D_d\cdot}_{p_d}$ and $\norm{\cdot}_{q_d}$ and maps $f_d(|x|)\eta$ to $c(x)\eta\,(1+|x|^{q_d})^{-(d-1)/2}$. Together with translations and scaling $\psi(x)\mapsto\sigma^{-(d-1)/2}\psi(x/\sigma)$ ($\sigma>0$), the Kelvin transform \eqref{eq:kelvin-spinor} generates the entire family \eqref{eq:optimizers} from $f_d(|x|)\eta$.

\subsection{Consequences}
On $\C\cong\R^2$ with identification $z=x_1+ix_2$ and Lebesgue measure $dA$, the Cauchy transform
\[
 (\mathcal Cg)(z)=\frac{1}{\pi}\int_\C\frac{g(\xi)}{z-\xi}\,dA(\xi)
\]
inverts $\dbar$ on $L^{4/3}(\C)$ and maps $L^{4/3}(\C)$ boundedly into $L^4(\C)$ \cite[Lemma~3.1(a)]{NRT}. As a corollary of the case $d=2$ of Theorem~\ref{thm}, we obtain its norm and all optimizers.

\begin{corollary}\label{cor:cauchy}
For every $g\in L^{4/3}(\C)$,
\begin{equation}\label{eq:cauchy-sharp}
 \norm{\mathcal Cg}_{L^4(\C)}\le\frac{2}{\pi}\norm{g}_{L^{4/3}(\C)},
\end{equation}
and the constant $2/\pi$ is sharp. Equality holds in \eqref{eq:cauchy-sharp} for $g\ne0$ if and only if, for some $\kappa\in\C\setminus\{0\}$,
\begin{equation}\label{eq:cauchy-eq}
 g(z)=\frac{\kappa}{\bigl(\sigma^4+|z-z_0|^4\bigr)^{3/2}}
 \qquad\text{or}\qquad
 g(z)=\frac{\kappa\,(z-z_1)|z-z_1|^2}{\bigl(|z-z_1|^4+|\alpha z+\beta|^4\bigr)^{3/2}},
\end{equation}
where $\sigma>0$ and $z_0\in\C$ in the first family, and $z_1\in\C$, $\alpha\ge0$, and $\beta\in\C$ with $\alpha z_1+\beta\ne0$ in the second.
\end{corollary}

Since $\dbar\,\mathcal Cg=g$, the functions in \eqref{eq:cauchy-eq} are the derivatives $\dbar u$ of planar optimizers $u$. At $\alpha=0$, the second family is a constant multiple of $\dbar\bigl((|\beta|^4+|z-z_1|^4)^{-1/2}\bigr)$, while the first family comes from its Kelvin image.

For a real-valued function $u$ on $\R^d$ and a unit spinor $\eta$, the spinor $\psi=c(\nabla u)\eta$ satisfies $|\psi|=|\nabla u|$ and $D_d\psi=-(\Delta u)\eta$. Theorem~\ref{thm} therefore implies a sharp $L^{p_d}\to L^{q_d}$ gradient estimate for the Laplacian.

\begin{corollary}\label{cor}
Let $d\in\{2,3,4\}$. For every $u\in L^1_{\mathrm{loc}}(\R^d;\R)$ with $\nabla u\in L^{q_d}(\R^d;\R^d)$ and $\Delta u\in L^{p_d}(\R^d)$,
\begin{equation}\label{eq:riesz-grad}
 \norm{\nabla u}_{L^{q_d}(\R^d)}\le S_d^{-1}\norm{\Delta u}_{L^{p_d}(\R^d)}.
\end{equation}
Equality holds for $\nabla u\not\equiv0$ if and only if
\begin{equation}\label{eq:riesz-eq}
 \nabla u(x)=\frac{C(x-x_0)}{\bigl(\sigma^{q_d}+|x-x_0|^{q_d}\bigr)^{(d-1)/2}}
\end{equation}
for some $C\in\R\setminus\{0\}$, $\sigma>0$, and $x_0\in\R^d$.
\end{corollary}

\subsection{The coupled Euler--Lagrange system}
Any critical point $\psi\in\X_d\setminus\{0\}$ of $Q_d$ satisfies the Euler--Lagrange equation
\begin{equation}\label{eq:EL}
 D_d\bigl(|D_d\psi|^{p_d-2}D_d\psi\bigr)=\Lambda\,|\psi|^{q_d-2}\psi,\qquad \Lambda=\frac{\norm{D_d\psi}_{p_d}^{p_d}}{\norm{\psi}_{q_d}^{q_d}}>0.
\end{equation}
With the normalization
\[
 U=\Lambda^{(d^2-1)/(4d)}\psi,\qquad V=\Lambda^{(d-1)^2/(4d)}|D_d\psi|^{p_d-2}D_d\psi,
\]
equation \eqref{eq:EL} becomes the first-order system
\begin{equation}\label{eq:coupled}
 D_dU=|V|^{q_d-2}V,\qquad D_dV=|U|^{q_d-2}U\qquad(U,V\in\X_d).
\end{equation}
Writing $\X_d^2=\X_d\times\X_d$, every solution $(U,V)\in\X_d^2$ of \eqref{eq:coupled} has common $L^{q_d}$ mass
\[
 I=\int_{\R^d}|U|^{q_d}\,dx=\int_{\R^d}|V|^{q_d}\,dx,
\]
with $I=Q_d(U)^d=Q_d(V)^d$ whenever $(U,V)\ne(0,0)$ (Lemma~\ref{lem:crit-norm}).

For $d\ge3$ and each unit spinor $\eta\in\C^{N_d}$, the trial spinor $f_d(|x|)\eta$ is itself a critical point: since $D_d(c(x)\eta)=-d\eta$, the pair
\begin{equation}\label{eq:trial}
 U_\eta(x)=d^{(d-1)/2}f_d(|x|)\,\eta,\qquad V_\eta(x)=-c(x)\,U_\eta(x)
\end{equation}
solves \eqref{eq:coupled} with $I=S_d^d$ (Lemma~\ref{lem:trial}). Let $\Or_d$ denote the orbit of $\{(U_\eta,V_\eta):|\eta|=1\}$ under translations, dilations, and the Kelvin transform $(U,V)\mapsto(U_{\mathcal K},-V_{\mathcal K})$. In dimensions three and four we bound $I$ from below for every nontrivial solution of the system and classify the solutions attaining the minimum.

\begin{theorem}\label{thm:pairs}
Let $d\in\{3,4\}$. Every nontrivial solution $(U,V)\in\X_d^2$ of \eqref{eq:coupled} satisfies $I\ge S_d^d$, with equality if and only if $(U,V)\in\Or_d$.
\end{theorem}

In particular, no critical point of $Q_d$ has quotient smaller than $S_d$ when $d\in\{3,4\}$, and the nontrivial solutions of \eqref{eq:coupled} with minimal $I$ are precisely the elements of $\Or_d$. Together with Proposition~\ref{prop:attain}, Theorem~\ref{thm:pairs} implies Theorem~\ref{thm} for $d=3,4$.

\subsection{Why two integral identities are needed}
In the plane the problem reduces to the scalar operator $2\dbar$. The Euler--Lagrange equation gives two complex-valued functions defining a map to $\mathbb{CP}^1$ that extends to an anti-holomorphic rational map of the sphere. By the area formula, the square of the best constant is $\pi^2$ times the rational degree in $\overline z$. The Aubin--Talenti trial gives the upper bound $\pi^2$, so the rational degree is one. The resulting M\"obius maps in $\overline z$ give all planar optimizers.

For the coupled system in dimensions three and four, we use a conformal metric built from both spinors. On the set where the pair is nonzero, put
\[
 \lambda=\frac{|U|^{q_d}+|V|^{q_d}}{|U|^2+|V|^2}.
\]
For the trial pair the factors involving $1+|x|^{q_d}$ cancel, giving $\lambda=d/(1+|x|^2)$ and hence the round metric $\lambda^2\,dx^2$. The Killing-spinor solution with $|U|=|V|$ gives the same metric. In dimension three the masses $I$ of these two pairs are $18\pi$ and $27\pi^2/4$, respectively. The volume of $\lambda^2\,dx^2$ therefore does not determine the mass.

The weighted Schr\"odinger--Lichnerowicz identity relates volume to scalar curvature and the twistor term. The Sobolev inequality gives a lower bound on volume when the remainder in this identity has nonnegative integral. Set $a=|U|^2/(|U|^2+|V|^2)$. With respect to $\lambda^d\,dx$, the density of $2I=\int_{\R^d}(|U|^{q_d}+|V|^{q_d})\,dx$ is $w(a)$. We choose a vector field whose weighted divergence contains $w(a)$ minus its mean on the trial pair. Integrating the divergence gives an identity for the difference between the normalized mass and volume.

In dimension three the inner product of the spinors contributes a term of unfavorable sign. Polarizing the weighted identity gives a current that controls this term, so the two integral identities can be used separately.

In dimension four we split the spinor space into its two half-spinor spaces. For a chiral pair the Clifford images of either nonzero spinor form a real basis of the opposite half-spinor space, so Bessel's inequality is an equality. The second remainder still has a negative coefficient near the ends of the amplitude range, so we estimate a linear combination of the two remainders.

When $I=S_d^d$, the quadratic defects vanish. The rescaled spinors then satisfy the Euclidean twistor equation and are affine. In dimension four the equality conditions for the chiral splitting recover the full spinor from its chiral parts. This yields the optimizer family in Theorem~\ref{thm}.

\section{Preliminaries}

We first establish the mapping properties of the inverses of $2\dbar$ and $D_d$. On $\C$ we write $\dot W^{1,4/3}(\C)=\{u\in L^4(\C):\nabla u\in L^{4/3}(\C)\}$, and define the integral operators
\[
 T_2g=\frac12\mathcal Cg,\qquad (T_dh)(x)=-\frac{1}{|\Sph^{d-1}|}\int_{\R^d}\frac{c(x-y)}{|x-y|^d}h(y)\,dy\quad(d\ge3).
\]

\begin{lemma}\label{lem:spaces}
Let $d\ge2$, and write $(T,\X,p,q)=(T_2,\dot W^{1,4/3}(\C),4/3,4)$ for $d=2$ and $(T_d,\X_d,p_d,q_d)$ for $d\ge3$. Then $T$ maps $L^p$ boundedly onto $\X$ as the two-sided inverse of $2\dbar$ or $D_d$, one has
\begin{equation}\label{eq:gradient-bound}
 \norm{\nabla Tg}_p+\norm{\nabla T^*g}_p\le C_d\norm{g}_p\qquad(g\in L^p),
\end{equation}
the space $\X_d$ coincides with $\dot W^{1,p_d}(\R^d;\C^{N_d})$ for $d\ge3$, and a function $u\in L^4(\C)$ has $\dbar u\in L^{4/3}(\C)$ distributionally if and only if $u\in\dot W^{1,4/3}(\C)$. Moreover, $C_c^\infty$ is dense in $\X$ under $\norm{\psi}_q+\norm{\nabla\psi}_p$.
\end{lemma}

\begin{proof}
For $d=2$, fractional integration \cite[Lemma~3.1(a)]{NRT} shows that $T_2$ maps $L^{4/3}(\C)$ boundedly into $L^4(\C)$, and from $\dbar(1/(\pi z))=\delta_0$ we have $(2\dbar)T_2=\mathrm{Id}$ on $L^{4/3}(\C)$. Its adjoint $T_2^*$ has kernel $-1/(2\pi\overline z)$ and inverts $-2\partial$. Since $\dbar T_2=\frac12\mathrm{Id}$ and $\partial T_2$ is half the Beurling transform, \eqref{eq:gradient-bound} follows from Calder\'on--Zygmund theory. If $u\in L^4(\C)$ has $\dbar u\in L^{4/3}(\C)$ distributionally, then $u-T_2(2\dbar u)$ is entire and belongs to $L^4(\C)$, hence vanishes by the mean-value property on disks; thus $u=T_2(2\dbar u)\in T_2(L^{4/3}(\C))=\dot W^{1,4/3}(\C)$.

For $d\ge3$, differentiating the Newtonian fundamental solution gives $D_dT_d=\mathrm{Id}$ on $C_c^\infty(\R^d)$ and, by continuity, distributionally on $L^{p_d}$. Since $|c(z)z'|=|z|\,|z'|$, the Hardy--Littlewood--Sobolev inequality \cite{Lieb} gives $T_d\colon L^{p_d}\to L^{q_d}$. The kernel also gives $T_d^*=T_d$ under real $L^{p_d}$--$L^{q_d}$ duality.

In Fourier variables, $D_d$ has symbol $ic(\xi)$, so $\partial_jT_d=\sum_{k=1}^d c_kR_jR_k$ in terms of the Riesz transforms $R_j$, whose $L^{p_d}$ boundedness proves \eqref{eq:gradient-bound}. If $\psi\in\X_d$, then $h=\psi-T_dD_d\psi\in L^{q_d}$ satisfies $D_dh=0$ and hence $\Delta h=0$; Weyl's lemma and the mean-value bound on balls imply $h=0$, so $\X_d=T_d(L^{p_d})=\dot W^{1,p_d}(\R^d;\C^{N_d})$.

Finally, multiplying by a radial cutoff $\chi(|\cdot|/R)$ and mollifying \cite[Section~4]{BM}, we see that $C_c^\infty$ is dense in $\X$ under $\norm{\psi}_q+\norm{\nabla\psi}_p$.
\end{proof}

By Lemma~\ref{lem:spaces}, $\norm{T_2}_{4/3\to4}^{-1}=\inf_{0\ne u\in\dot W^{1,4/3}(\C)}\norm{2\dbar u}_{4/3}/\norm{u}_4$ and $\norm{T_d}_{p_d\to q_d}^{-1}=\inf_{\X_d\setminus\{0\}}Q_d$ are finite and positive. We use Ledoux's improved Sobolev inequality \cite{Ledoux} and a duality argument to prove attainment. For the final norm splitting, compare \cite[Lemma~2.7]{Lieb}. See \cite{Lions} for the concentration-compactness approach.

\begin{proposition}\label{prop:attain}
The operators $T_2\colon L^{4/3}(\C)\to L^4(\C)$ and $T_d\colon L^{p_d}\to L^{q_d}$ for $d\ge3$ attain their norms. Consequently, $\inf_{0\ne u\in\dot W^{1,4/3}(\C)}\norm{2\dbar u}_{4/3}/\norm{u}_4$ is attained, as is $\inf_{\X_d\setminus\{0\}}Q_d$ for every $d\ge2$, using Lemma~\ref{lem:d2-scalar} for $d=2$.
\end{proposition}

\begin{proof}
Write $(T,\X,p,q)$ as in Lemma~\ref{lem:spaces} and put $M_T=\norm{T}_{p\to q}>0$. By Ledoux's inequality \cite[Theorem~1]{Ledoux} applied to each real component and density of $C_c^\infty$ in $\X$, we have with $\theta=p/q=(d-1)/(d+1)$
\begin{equation}\label{eq:ledoux}
 \norm{u}_q\le C\norm{\nabla u}_p^{\theta}\Bigl(\sup_{\tau>0}\tau^{(d-1)/4}\norm{e^{\tau\Delta}u}_\infty\Bigr)^{1-\theta}\qquad(u\in\X).
\end{equation}
Let $(g_n)$ satisfy $\norm{g_n}_p=1$ and $\norm{u_n}_q\to M_T$, where $u_n=Tg_n$. By \eqref{eq:gradient-bound} and \eqref{eq:ledoux}, there are $x_n$ and $\tau_n>0$ with $\tau_n^{(d-1)/4}|e^{\tau_n\Delta}u_n(x_n)|\ge c_0>0$. Rescaling $(g_n,u_n)$ by $r_n=\sqrt{\tau_n}$ around $x_n$,
\[
 \bigl(g_n,u_n\bigr)\longmapsto\Bigl(r_n^{(d+1)/2}g_n(x_n+r_n\cdot),\,r_n^{(d-1)/2}u_n(x_n+r_n\cdot)\Bigr),
\]
we may assume that $|e^{\Delta}u_n(0)|\ge c_0$. Along a subsequence, $g_n\rightharpoonup g$ in $L^p$ and $u_n\rightharpoonup u=Tg$ in $L^q$. Pairing $u_n$ with the heat kernel at $\tau=1$ shows that $u\ne0$, so $0<\norm{g}_p\le1$, and by the local Rellich compactness theorem we may also assume that $u_n\to u$ almost everywhere.

We next show that $g_n$ converges almost everywhere along a subsequence. For $z\in L^q\setminus\{0\}$, set $J_q(z)=|z|^{q-2}z/\norm{z}_q^{q-1}$, so that $\norm{J_q(z)}_p=1$ and $\int\ip{z}{J_q(z)}\,dx=\norm{z}_q$. For large $n$, put $h_n=T^*J_q(u_n)$. Since $\norm{T^*}_{p\to q}=M_T$, H\"older's inequality gives
\[
 \norm{u_n}_q=\int\ip{g_n}{h_n}\,dx\le\norm{h_n}_q\le M_T,
\]
and therefore $\norm{h_n}_q\to M_T$. With $\tilde g_n=J_q(h_n)$, we have
\[
 \frac{\norm{u_n}_q+\norm{h_n}_q}{2\norm{h_n}_q}
 =\int\Bigl\langle\frac{g_n+\tilde g_n}{2},\frac{h_n}{\norm{h_n}_q}\Bigr\rangle\,dx
 \le\biggl\lVert\frac{g_n+\tilde g_n}{2}\biggr\rVert_p\le1,
\]
so $\norm{(g_n+\tilde g_n)/2}_p\to1$. The uniform convexity of $L^p$ then gives $\norm{g_n-\tilde g_n}_p\to0$.

By \eqref{eq:gradient-bound}, $h_n$ is bounded in $W^{1,p}_{\mathrm{loc}}$. Passing to a subsequence, we may assume that $h_n\to h$ almost everywhere. Since $\norm{h_n}_q\to M_T>0$, we have $\tilde g_n\to |h|^{q-2}h/M_T^{q-1}$ almost everywhere. Passing to a further subsequence in $g_n-\tilde g_n\to0$ in $L^p$ gives the same almost-everywhere limit for $g_n$. Boundedness in $L^p$ identifies this limit with the weak limit $g$.

We now apply the Brezis--Lieb lemma \cite[Theorem~1]{BL} to $g_n$ and $u_n$. For finite-dimensional spinors, the proof uses the same estimate
\[
 \bigl||v+w|^r-|v|^r\bigr|\le\varepsilon |v|^r+C_{r,\varepsilon}|w|^r
 \qquad(r>1),
\]
which follows from the triangle inequality and the mean value theorem, so the lemma applies to the spinor norm. By weak convergence, $\norm{g}_p\le1$, and hence
\[
 M_T^q\le M_T^q\bigl(\norm{g}_p^q+(1-\norm{g}_p^p)^{q/p}\bigr).
\]
Since $q/p>1$ and $0<\norm{g}_p^p\le1$, we conclude that $\norm{g}_p=1$ and $\norm{Tg}_q=M_T$.
\end{proof}

When $d\ge2$ is fixed, we abbreviate $p=p_d$, $q=q_d$, and $\nu=q_d-2=2/(d-1)$ (with $|z'|^{p-2}z'=0$ at $z'=0$).

\begin{lemma}\label{lem:crit-norm}
Let $d\ge2$.
\begin{enumerate}
\item[\textup{(i)}] Every critical point $\psi\in\X_d\setminus\{0\}$ of $Q_d$ solves \eqref{eq:EL}, and the rescaled spinors
\begin{equation}\label{eq:normalize}
 U=\Lambda^{(d^2-1)/(4d)}\psi,\qquad V=\Lambda^{(d-1)^2/(4d)}|D_d\psi|^{p-2}D_d\psi
\end{equation}
solve \eqref{eq:coupled} in $\X_d^2$ with $\norm{U}_q^q=\norm{V}_q^q=Q_d(\psi)^d$.
\item[\textup{(ii)}] Conversely, every solution $(U,V)\in\X_d^2$ of \eqref{eq:coupled} satisfies
\begin{equation}\label{eq:equal-norm}
 I=\int_{\R^d}|U|^q\,dx=\int_{\R^d}|V|^q\,dx,
\end{equation}
with $I=Q_d(U)^d=Q_d(V)^d$ whenever $(U,V)\ne(0,0)$. If one component vanishes identically, so does the other.
\end{enumerate}
\end{lemma}

\begin{proof}
\textup{(i)} Since $z'\mapsto|z'|^p/p$ belongs to $C^1(\C^{N_d})$ with differential $\ip{|z'|^{p-2}z'}{\cdot}$, real variations of $Q_d$ at $\psi$ lead to \eqref{eq:EL} \cite[(5.2), (C.1)]{WZ}, and direct scaling shows that $(U,V)$ in \eqref{eq:normalize} solves \eqref{eq:coupled} with $\norm{U}_q^q=\norm{V}_q^q=Q_d(\psi)^d$.

\textup{(ii)} Let $(U,V)\in\X_d^2$ solve \eqref{eq:coupled}, and set $j_k=\ip{c_kV}{U}\in L^{d/(d-1)}(\R^d)$. Choose $\chi\in C_c^\infty(\R)$ equal to $1$ on $[-1,1]$ and $0$ outside $(-2,2)$, and put $\chi_R(x)=\chi(|x|/R)$. Approximating $U$ and $V$ in $\X_d$ by $C_c^\infty$ spinors (Lemma~\ref{lem:spaces}) and using $c_k^*=-c_k$, we find
\[
 -\int_{\R^d}\nabla\chi_R\cdot j\,dx=\int_{\R^d}\chi_R\bigl(\ip{D_dV}{U}-\ip{V}{D_dU}\bigr)\,dx=\int_{\R^d}\chi_R\bigl(|U|^q-|V|^q\bigr)\,dx.
\]
On $A_R=\{R<|x|<2R\}$, H\"older's inequality bounds the left-hand side by $C\norm{j}_{L^{d/(d-1)}(A_R)}$, which tends to $0$ as $R\to\infty$. This proves \eqref{eq:equal-norm}, and the remaining assertions follow from $|D_dU|^p=|V|^q$ and $|D_dV|^p=|U|^q$.
\end{proof}

\begin{lemma}\label{lem:removable}
Let $d\ge2$, $x_0\in\R^d$, and $L=2\dbar$ for $d=2$ or $L=D_d$ for $d\ge3$. If $W\in L^{q_d}_{\mathrm{loc}}(\R^d)$ satisfies $LW=f$ distributionally on $\R^d\setminus\{x_0\}$ with $f\in L^{p_d}_{\mathrm{loc}}(\R^d)$, then $LW=f$ distributionally on $\R^d$.
\end{lemma}

\begin{proof}
Testing against $\bigl(1-\chi(|x-x_0|/\eps)\bigr)\varphi$ for $\varphi\in C_c^\infty(\R^d)$, with $\chi$ as above, and using $d/p_d-1=(d-1)/2>0$, we see by H\"older's inequality that the cutoff commutator is $O\bigl(\eps^{(d-1)/2}\norm{W}_{L^{q_d}(B_{2\eps}(x_0))}\bigr)$, which tends to $0$ as $\eps\downarrow0$ (cf.\ \cite[(49)]{BF}).
\end{proof}

\begin{lemma}\label{lem:trial}
Let $d\ge3$. For every unit spinor $\eta\in\C^{N_d}$, the pair $(U_\eta,V_\eta)$ in \eqref{eq:trial} belongs to $\X_d^2$ and solves \eqref{eq:coupled} with $I=S_d^d$. In particular, $\inf Q_d\le S_d$.
\end{lemma}

\begin{proof}
Using $f_d'(r)=-dr^{q-1}(1+r^q)^{-(d+1)/2}$ and $D_d(c(x)\eta)=-d\eta$, we compute on $\R^d\setminus\{0\}$
\[
 D_dU_\eta=d^{(d-1)/2}\,\frac{f_d'(r)}{r}\,c(x)\eta=|V_\eta|^{q-2}V_\eta,\qquad
 D_dV_\eta=d^{(d-1)/2}\bigl(df_d(r)+rf_d'(r)\bigr)\eta=|U_\eta|^{q-2}U_\eta.
\]
Since $U_\eta,V_\eta\in L^q(\R^d)$, both equations extend across $0$ by Lemma~\ref{lem:removable}. In polar coordinates with $s=r^q$, \eqref{eq:Sd} implies $\int_{\R^d}|U_\eta|^q\,dx=S_d^d$ \cite[Theorem~4.1]{BM}, and $Q_d(U_\eta)=S_d$ follows from \eqref{eq:equal-norm}.
\end{proof}

Finally, we record the conformal symmetries of \eqref{eq:coupled}. Recall that $\Or_d$ is the orbit of $\{(U_\eta,V_\eta):|\eta|=1\}$ under translations, dilations $W(x)\mapsto\sigma^{-(d-1)/2}W(x/\sigma)$ ($\sigma>0$), and the Kelvin transform
\begin{equation}\label{eq:kelvin}
 (U,V)\longmapsto(U_{\mathcal K},-V_{\mathcal K}).
\end{equation}

\begin{lemma}\label{lem:conf}
Let $d\ge3$. Translations, dilations, and \eqref{eq:kelvin} map solutions of \eqref{eq:coupled} in $\X_d^2$ to solutions of \eqref{eq:coupled} in $\X_d^2$, preserve $\norm{\cdot}_q$ and $\norm{D_d\cdot}_p$ of each component, and preserve $I$ and the orbit $\Or_d$.
\end{lemma}

\begin{proof}
Translations and dilations are immediate by scaling. For $x\in\R^d\setminus\{0\}$, write $n=x/|x|$. By the pointwise spinorial Kelvin identity \cite[Lemma~4.1]{BF} (cf.\ \cite[Remark~3]{FL}) with $c_j^*=-c_j$,
\begin{equation}\label{eq:kelvin-cov}
 D_d(W_{\mathcal K})(x)=-|x|^{-d-1}c(n)\,(D_dW)\bigl(x/|x|^2\bigr)\qquad(x\ne0),
\end{equation}
and the Jacobian $|x|^{-2d}$ preserves $\norm{\cdot}_q$. Since $|W_{\mathcal K}(x)|=|x|^{1-d}|W(x/|x|^2)|$, we also have $|W_{\mathcal K}|^{q-2}W_{\mathcal K}(x)=|x|^{-d-1}c(n)(|W|^{q-2}W)(x/|x|^2)$ off the origin. Together with \eqref{eq:kelvin-cov} and Lemma~\ref{lem:removable} at $x_0=0$, this shows that $(U_{\mathcal K},-V_{\mathcal K})\in\X_d^2$ solves \eqref{eq:coupled} on $\R^d$ with the same norms and the same value of $I$.
\end{proof}

\section{The planar inequality}\label{sec:d2}

Throughout this section $d=2$, $z=x_1+ix_2$, and $p=4/3$. For any skew-Hermitian Clifford pair $(c_1,c_2)$ on $\C^2$, the chirality operator $\Gamma=ic_1c_2$ satisfies $\Gamma^*=\Gamma$, $\Gamma^2=I_2$, and $c_1\Gamma=-\Gamma c_1$. In the orthonormal basis $(e_+,ic_1e_+)$ built from a unit $+1$ eigenvector $e_+$ of $\Gamma$, we have
\begin{equation}\label{eq:d2-rep}
 c(x)=-i\begin{pmatrix}0&\overline z\\ z&0\end{pmatrix},\qquad
 D_2=-2i\begin{pmatrix}0&\partial\\ \dbar&0\end{pmatrix}.
\end{equation}
A unitary basis change preserves the norms and conjugates $c(x)$, so we may work in the representation \eqref{eq:d2-rep}. Set
\[
 S_*=\inf_{0\ne u\in\dot W^{1,4/3}(\C)}\frac{\norm{2\dbar u}_{4/3}}{\norm{u}_4}.
\]

\begin{lemma}\label{lem:d2-scalar}
One has $\inf_{\X_2\setminus\{0\}}Q_2=S_*\le\pi=S_2$. For $\psi=(u,v)\in\X_2\setminus\{0\}$ with $u\ne0$ and $v\ne0$, equality $Q_2(\psi)=S_*$ holds if and only if both $u$ and $\overline v$ attain $S_*$ and $|u|/|v|$ and $|\dbar u|/|\dbar\overline v|$ are constant almost everywhere.
\end{lemma}

\begin{proof}
For $\psi=(u,v)\in\X_2$, Lemma~\ref{lem:spaces} shows that both $u$ and $\overline v$ belong to $\dot W^{1,4/3}(\C)$, with $|D_2\psi|^2=|2\dbar u|^2+|2\dbar\overline v|^2$. By the reverse Minkowski inequality in $L^{2/3}$ and Minkowski's inequality in $L^2$,
\[
 \norm{D_2\psi}_{4/3}^2\ge\norm{2\dbar u}_{4/3}^2+\norm{2\dbar\overline v}_{4/3}^2\ge S_*^2\bigl(\norm{u}_4^2+\norm{v}_4^2\bigr)\ge S_*^2\norm{\psi}_4^2.
\]
When $u,v\ne0$, strict concavity of $s\mapsto s^{2/3}$ and strict convexity of $s\mapsto s^2$ show that equality holds in the outer steps if and only if $|2\dbar u|$ is proportional to $|2\dbar\overline v|$ and $|u|$ is proportional to $|v|$ almost everywhere. Taking $\psi=(u,0)$ shows that $\inf Q_2=S_*$. Finally, for $f_2(z)=(1+|z|^4)^{-1/2}$ we have $\dbar f_2=-z|z|^2f_2^3$ and, in polar coordinates with $s=r^2$ and $r\mapsto1/r$,
\begin{equation}\label{eq:planar-J}
 \int_\C|f_2|^4\,dA=\pi\int_0^\infty\frac{ds}{(1+s^2)^2}=\frac{\pi^2}{4}=\int_\C|z|^4|f_2|^4\,dA.
\end{equation}
Thus $\norm{2\dbar f_2}_{4/3}=2(\pi^2/4)^{3/4}=\pi\norm{f_2}_4$, and so $S_*\le\pi=S_2$.
\end{proof}

\begin{lemma}\label{lem:d2-cubic}
Let $u\in\dot W^{1,4/3}(\C)$ satisfy $\norm{u}_4=1$ and $\norm{2\dbar u}_{4/3}=S_*$, and put $W=|2\dbar u|^{-2/3}(2\dbar u)$ \textup{(}with $W=0$ where $\dbar u=0$\textup{)}. Then $(P,Q)=\bigl(2^{-1/2}S_*^{1/2}u,\,-2^{-1/2}S_*^{1/6}\overline W\bigr)$ belongs to $C^\infty(\C;\C^2)$ and solves
\begin{equation}\label{eq:cubic}
 \dbar P=-Q\overline Q^{\,2},\qquad \dbar Q=P\overline P^{\,2},\qquad \int_\C|P|^4\,dA=\int_\C|Q|^4\,dA=\frac{S_*^2}{4}.
\end{equation}
\end{lemma}

\begin{proof}
By definition, $2\dbar u=|W|^2W$ and $\int_\C|W|^4\,dA=\int_\C|2\dbar u|^{4/3}\,dA=S_*^{4/3}$. For any $\varphi\in C_c^\infty(\C;\C)$, differentiating $\norm{2\dbar(u+\tau\varphi)}_{4/3}/\norm{u+\tau\varphi}_4$ at $\tau=0$ along both $\varphi$ and $i\varphi$, we find
\[
 \int_\C\overline W\,(2\dbar\varphi)\,dA=S_*^{4/3}\int_\C|u|^2\overline u\,\varphi\,dA.
\]
Hence $2\dbar\overline W=-S_*^{4/3}|u|^2\overline u$ distributionally on $\C$, and $(P,Q)$ solves \eqref{eq:cubic}. In the representation \eqref{eq:d2-rep}, $\Psi=(P,i\overline Q)\in L^4(\R^2;\C^2)$ solves $D_2\Psi=h(\Psi)\Psi$ with $h(z)=2\operatorname{diag}(|z_1|^2,|z_2|^2)$. This is the nonlinearity in \cite[(50)]{BF} with $\beta_1=2$ and $\beta_2=0$, so \cite[Proposition~5.3]{BF} gives $P,Q\in C^\infty(\C)$.
\end{proof}

On the pair coming from $f_2$, the ratio $Q/P$ equals $\overline z$. Clearing the denominator in $\partial(Q/P)$ leads to the Wronskian $P\partial Q-Q\partial P$.

\begin{lemma}\label{lem:planar-proj}
For $(P,Q)$ as in Lemma~\ref{lem:d2-cubic}, the Wronskian $P\partial Q-Q\partial P$ vanishes identically on $\C$, and $[P:Q]\colon\C\setminus\{P=Q=0\}\to\mathbb{CP}^1$ extends to a nonconstant anti-holomorphic rational map $\widehat\C\to\mathbb{CP}^1$ of rational degree one in $\overline z$. In particular, $S_*=\pi$.
\end{lemma}

\begin{proof}
Let $u$ be any normalized minimizer for $S_*$, which exists by Proposition~\ref{prop:attain}, and let $(P,Q)$ be as in Lemma~\ref{lem:d2-cubic}. By Lemma~\ref{lem:spaces}, $\partial P,\partial Q\in L^{4/3}(\C)$, so $P\partial Q-Q\partial P\in L^1(\C)$ by H\"older's inequality. Using $\partial\overline P=-Q^2\overline Q$ and $\partial\overline Q=P^2\overline P$ from \eqref{eq:cubic}, we compute
\[
 \dbar(P\partial Q-Q\partial P)=2P^2\overline P\,\partial\overline P+2Q^2\overline Q\,\partial\overline Q=0.
\]
Being entire and integrable on $\C$, $P\partial Q-Q\partial P$ vanishes identically. Wherever $P\ne0$, the quotient $\zeta=Q/P$ therefore satisfies $\partial\zeta=0$ and
\begin{equation}\label{eq:planar-zeta}
 \dbar\zeta=\overline P^{\,2}\bigl(1+|\zeta|^4\bigr).
\end{equation}
Similarly $\partial(P/Q)=0$ wherever $Q\ne0$, so $[P:Q]$ is anti-holomorphic off common zeros of $P$ and $Q$.

To extend $[P:Q]$ across common zeros (cf.\ \cite[Appendix, Lemma~A.1]{BC}), write $F=(P,Q)^T$ and set
\[
 M_0=\frac{(\dbar F)F^*}{|F|^2}\quad\text{on }\{F\ne0\},\qquad M_0=0\quad\text{on }\{F=0\},
\]
so that $\dbar F=M_0F$ and $|M_0|\le C|F|^2\in L^\infty_{\mathrm{loc}}(\C)$. Choose $\delta>0$ on a disk $B=B_\delta(z_0)$ so that $4\delta\norm{M_0}_{L^\infty(B)}<1/4$. With $\chi_B$ the indicator of $B$, the map $G\mapsto I_2+\mathcal C(M_0\chi_BG)$ is a contraction on $L^\infty(B;\C^{2\times2})$, and its continuous fixed point $G$ satisfies $\norm{G-I_2}_{L^\infty(B)}<1/2$.

The Cauchy transform of the bounded, compactly supported $M_0\chi_BG$ is continuous and locally $W^{1,r}$ for every finite $r>2$. The product rule therefore gives $\dbar(G^{-1}F)=0$ on $B$. If $F$ vanished on an open set, we could choose such a disk $B=B_\delta(z_0)$ at a boundary point of its interior zero set. The identity theorem would then give $G^{-1}F=0$ on $B$, contradicting the choice of its center. Since $F\not\equiv0$, the holomorphic pair $G^{-1}F$ does not vanish identically on $B$. The common zeros of $P$ and $Q$ are therefore isolated.

At a common zero, factor out the minimal vanishing power $(z-z_0)^m$ of the two components of $G^{-1}F$. Then $[F]=[G(G^{-1}F)]$ has a limit at $z_0$, and Riemann's removable singularity theorem extends $[P:Q]$ anti-holomorphically across $z_0$.

Near infinity, set $F_{\mathcal K}(z)=(i/z)F(1/z)$ on $\C\setminus\{0\}$. Then $F_{\mathcal K}\in L^4(\C;\C^2)$ solves \eqref{eq:cubic} off $0$, hence on $\C$ by Lemma~\ref{lem:removable}, and is smooth by the regularity argument in Lemma~\ref{lem:d2-cubic}. Applying the preceding argument to $F_{\mathcal K}$ at $0$ extends $[F(1/z)]=[F_{\mathcal K}(z)]$ across $0$, and thus $[P:Q]$ across $\infty$. If $[P:Q]$ were constant, then \eqref{eq:planar-zeta} (or its analogue for $P/Q$) would imply $P\equiv0$ or $Q\equiv0$, contradicting $\int_\C|P|^4\,dA=\int_\C|Q|^4\,dA=S_*^2/4>0$.

It follows that $\zeta(z)$ is a rational function of $\overline z$ of degree $k\ge1$. Equation \eqref{eq:planar-zeta} gives $|P|^4=|\dbar\zeta|^2(1+|\zeta|^4)^{-2}$, with $|\dbar\zeta|^2$ the absolute Jacobian of $\zeta$. The area formula and \eqref{eq:planar-J} yield
\[
 \frac{S_*^2}{4}=\int_\C|P|^4\,dA=k\int_\C\frac{dA(\xi)}{(1+|\xi|^4)^2}=\frac{k\pi^2}{4}.
\]
Since $S_*\le\pi$ by Lemma~\ref{lem:d2-scalar}, we conclude that $k=1$ and $S_*=\pi$.
\end{proof}

\begin{proposition}\label{prop:d2-eq}
A nonzero function $u\in\dot W^{1,4/3}(\C)$ attains $\norm{2\dbar u}_{4/3}/\norm{u}_4=\pi$ if and only if
\begin{equation}\label{eq:planar-family}
 u(z)=\frac{\kappa\,\ell_1(z)}{\sqrt{\mathcal L(z)}},\qquad \mathcal L=|\ell_1|^4+|\ell_2|^4,
\end{equation}
for some $\kappa\in\C\setminus\{0\}$ and $\ell_j(z)=\alpha_j\overline z+\beta_j$ \textup{($j=1,2$)} with $\mathcal W=\alpha_2\beta_1-\alpha_1\beta_2\ne0$. Consequently, in the representation \eqref{eq:d2-rep}, $\psi\in\X_2\setminus\{0\}$ attains $Q_2(\psi)=\pi$ if and only if
\begin{equation}\label{eq:d2-spinor-opt}
 \psi(z)=\frac{1}{\sqrt{\mathcal L(z)}}\bigl(\kappa_1\ell_1(z),\,\kappa_2\overline{\ell_1(z)}\bigr)
\end{equation}
with $(\kappa_1,\kappa_2)\in\C^2\setminus\{0\}$ and $\ell_1,\ell_2$ as above.
\end{proposition}

\begin{proof}
Suppose $u\in\dot W^{1,4/3}(\C)\setminus\{0\}$ attains $S_*=\pi$, normalized so that $\norm{u}_4=1$. By Lemmas~\ref{lem:d2-cubic} and~\ref{lem:planar-proj}, the associated pair $(P,Q)$ solves \eqref{eq:cubic} with $k=1$, so $\zeta=Q/P=\ell_2/\ell_1$ with $\mathcal W=\alpha_2\beta_1-\alpha_1\beta_2\ne0$ and $\mathcal L=|\ell_1|^4+|\ell_2|^4>0$ on $\C$. Inserting $\dbar\zeta=\mathcal W/\ell_1^2$ and $1+|\zeta|^4=\mathcal L/|\ell_1|^4$ into \eqref{eq:planar-zeta}, we obtain
\[
 P^2=\frac{\overline{\mathcal W}\,\ell_1^2}{\mathcal L}
\]
wherever $P\ne0$ and $\ell_1\ne0$. On the connected open set $\{\ell_1\ne0\}$, $P$ cannot vanish: at any zero $z_0\in\{\ell_1\ne0\}$, continuity of $Q=(\ell_2/\ell_1)P$ would make $z_0$ a common zero of $(P,Q)$, whereas $\lim_{z\to z_0}P^2=\overline{\mathcal W}\,\ell_1(z_0)^2/\mathcal L(z_0)\ne0$. Thus $P\sqrt{\mathcal L}/\ell_1$ is constant on $\{\ell_1\ne0\}$, and \eqref{eq:planar-family} follows by continuity.

Conversely, for $u$ as in \eqref{eq:planar-family}, using $\dbar\mathcal L=2\sum_{j=1}^2\alpha_j\ell_j\overline{\ell_j}^{\,2}$ we compute
\begin{equation}\label{eq:planar-A}
 2\dbar\Bigl(\frac{\kappa\ell_1}{\sqrt{\mathcal L}}\Bigr)=-\frac{2\kappa\mathcal W\ell_2\overline{\ell_2}^{\,2}}{\mathcal L^{3/2}}.
\end{equation}
The map $\zeta=\ell_2/\ell_1$ is one-to-one off $\{\ell_1=0\}$, and its image omits at most one point of $\C$. Changing variables with absolute Jacobian $|\mathcal W|^2/|\ell_1|^4$ and using \eqref{eq:planar-J}, we obtain
\[
 \int_\C|u|^4\,dA=\frac{|\kappa|^4}{|\mathcal W|^2}\,\frac{\pi^2}{4},\qquad \int_\C|2\dbar u|^{4/3}\,dA=\frac{(2|\kappa\mathcal W|)^{4/3}}{|\mathcal W|^2}\,\frac{\pi^2}{4},
\]
so $\norm{2\dbar u}_{4/3}/\norm{u}_4=\pi$.

Finally, let $\psi=(u,v)\in\X_2\setminus\{0\}$ attain $Q_2(\psi)=\pi$. If $u=0$ or $v=0$, the remaining component (or its conjugate) has the form \eqref{eq:planar-family}. If both $u,v\ne0$, Lemma~\ref{lem:d2-scalar} shows that $u=\kappa_1'\ell_1/\sqrt{\mathcal L}$ and $\overline v=\kappa_2'\widetilde\ell_1/\sqrt{\widetilde{\mathcal L}}$ with $|u|/|v|$ constant. Comparing zero sets in $\C$, we see that $\widetilde\ell_1$ is a nonzero complex multiple of $\ell_1$; hence $\widetilde{\mathcal L}$ is a constant multiple of $\mathcal L$ and $\psi$ has the form \eqref{eq:d2-spinor-opt}. The converse is immediate since $|\psi|=(|\kappa_1|^2+|\kappa_2|^2)^{1/2}|\ell_1|/\sqrt{\mathcal L}$ and $|D_2\psi|=(|\kappa_1|^2+|\kappa_2|^2)^{1/2}|2\dbar(\ell_1/\sqrt{\mathcal L})|$.
\end{proof}

\begin{proof}[Proof of Theorem~\ref{thm} for $d=2$]
By Lemmas~\ref{lem:d2-scalar} and~\ref{lem:planar-proj} and Proposition~\ref{prop:d2-eq}, we only need to match \eqref{eq:d2-spinor-opt} with \eqref{eq:optimizers} at $d=2$. When $\alpha_1=0$, we have $\ell_1=\beta_1\ne0$ and $|\ell_2|=|\alpha_2|\,|z-z_0|$ with $z_0=-\overline{\beta_2}/\overline{\alpha_2}$; this is the first family in \eqref{eq:optimizers} with $\sigma=|\beta_1/\alpha_2|>0$ and $\eta=|\alpha_2|^{-2}(\kappa_1\beta_1,\kappa_2\overline{\beta_1})$. When $\alpha_1\ne0$, we write $\ell_1=\alpha_1(\overline z-\overline{z_1})$ with $z_1=-\overline{\beta_1}/\overline{\alpha_1}$; comparing with $c(x-x_1)$ in \eqref{eq:d2-rep}, we obtain the second family in \eqref{eq:optimizers} with $\alpha=|\alpha_2/\alpha_1|$, $|\alpha x+\beta|=|\ell_2/\alpha_1|$ (so $|\alpha z_1+\beta|=|\mathcal W|/|\alpha_1|^2\ne0$), and $\eta=i|\alpha_1|^{-2}(\kappa_2\overline{\alpha_1},\kappa_1\alpha_1)$. Conversely, taking $(\ell_1,\ell_2)=(\sigma,\overline z-\overline{z_0})$ in the first family and $(\ell_1,\ell_2)=(\overline z-\overline{z_1},\,\alpha\overline z+\overline\beta)$ in the second realizes every spinor in \eqref{eq:optimizers} as \eqref{eq:d2-spinor-opt} with $\mathcal W\ne0$.
\end{proof}

\begin{proof}[Proof of Corollary~\ref{cor:cauchy}]
In the representation \eqref{eq:d2-rep}, for any $g\in L^{4/3}(\C)$ the spinor $\psi=(\mathcal Cg,0)\in\X_2$ satisfies $\dbar\,\mathcal Cg=g$ and $|D_2\psi|=2|g|$. By Proposition~\ref{prop:d2-eq}, $\norm{\mathcal Cg}_4\le(2/\pi)\norm{g}_{4/3}$, with equality for $g\ne0$ if and only if $\mathcal Cg=\kappa_1\ell_1/\sqrt{\mathcal L}$ as in \eqref{eq:planar-family}. Since $u=T_2(2\dbar u)=\mathcal C(\dbar u)$ on $\dot W^{1,4/3}(\C)$ (Lemma~\ref{lem:spaces}), equality holds if and only if, by \eqref{eq:planar-A},
\[
 g=\dbar\Bigl(\frac{\kappa_1\ell_1}{\sqrt{\mathcal L}}\Bigr)=-\frac{\kappa_1\mathcal W\,\overline{\ell_2}|\ell_2|^2}{\mathcal L^{3/2}}.
\]
When $\alpha_2=0$, we have $\beta_2\ne0$, $\alpha_1\ne0$, and $|\ell_1|=|\alpha_1|\,|z-z_0|$ with $z_0=-\overline{\beta_1}/\overline{\alpha_1}$, which is the first family in \eqref{eq:cauchy-eq} with $\sigma=|\beta_2/\alpha_1|>0$. When $\alpha_2\ne0$, writing $\overline{\ell_2}=\overline{\alpha_2}(z-z_1)$ with $z_1=-\overline{\beta_2}/\overline{\alpha_2}$ and $|\ell_1/\alpha_2|=|\alpha z+\beta|$ with $\alpha=|\alpha_1/\alpha_2|\ge0$ and $|\alpha z_1+\beta|=|\mathcal W|/|\alpha_2|^2\ne0$, we obtain the second family in \eqref{eq:cauchy-eq}. Conversely, setting $(\ell_1,\ell_2)=(\overline z-\overline{z_0},\sigma)$ in the first family and $(\ell_1,\ell_2)=(\alpha\overline z+\overline\beta,\,\overline z-\overline{z_1})$ in the second realizes every function in \eqref{eq:cauchy-eq}.
\end{proof}

\section{Integral identities}\label{sec:euclidean}

Throughout this section we fix a dimension $d\ge3$, write $p=p_d$, $q=q_d$, and $\nu=q_d-2=2/(d-1)$, and let $(U,V)\in\X_d^2$ be a nontrivial solution of \eqref{eq:coupled} on $\R^d$. We derive the integral identities \eqref{eq:vol-id} and \eqref{eq:second-id}. Since $\lambda$ vanishes where $U=V=0$, we first establish the regularity and weighted gradient bound needed to integrate the divergences.

\subsection{The weighted identity}
For any weakly differentiable spinor $\varphi$ on an open subset of $\R^d$, the twistor operator is defined by
\[
 \T_j\varphi=\partial_j\varphi+\frac1dc_jD_d\varphi\qquad(j=1,\dots,d).
\]
From the Clifford relations, $\sum_{j=1}^d c_j\T_j\varphi=0$ and
\begin{equation}\label{eq:twistor-norm}
 |\T\varphi|^2=\sum_{j=1}^d|\T_j\varphi|^2=|\nabla\varphi|^2-\frac1d|D_d\varphi|^2.
\end{equation}

We first prove a weighted Schr\"odinger--Lichnerowicz identity. We apply it with a regularized weight to control derivatives near the zero set of $(U,V)$, and then with the conformal factor $\lambda$.

\begin{lemma}\label{lem:sl}
Let $d\ge3$ and let $O\subset\R^d$ be open. Let $\Phi\in C^2(O;\C^{N_d})$ and let $u\in C^2(O)$ be positive. Write $\zeta=u^{(d-2)/2}$ and $\ip{c\,W}{Z}=(\ip{c_jW}{Z})_{j=1}^d$. Then on $O$,
\begin{equation}\label{eq:sl}
 \begin{aligned}
  u^{d-1}\bigl|\T(u^{-d/2}\Phi)\bigr|^2-\frac{d-1}{d}u^{-1}|D_d\Phi|^2
  &=\frac{d-1}{d-2}\bigl(u^{1-d}|\Phi|^2\bigr)\zeta\Delta\zeta\\
  &\quad+\Div\Bigl(u^{-1}\ip{c\,D_d\Phi}{\Phi}+\frac12u^{d-2}\nabla\bigl(u^{1-d}|\Phi|^2\bigr)\Bigr).
 \end{aligned}
\end{equation}
\end{lemma}

\begin{proof}
With $\varphi=u^{-d/2}\Phi$, we have
\[
 \partial_j\varphi=u^{-d/2}\Bigl(\partial_j\Phi-\frac{d}{2}u^{-1}(\partial_ju)\Phi\Bigr),\qquad
 D_d\varphi=u^{-d/2}\Bigl(D_d\Phi-\frac{d}{2}u^{-1}c(\nabla u)\Phi\Bigr).
\]
Inserting both expressions into \eqref{eq:twistor-norm} and using $c_j^*=-c_j$, we find
\begin{align*}
 u^{d-1}|\T\varphi|^2-\frac{d-1}{d}u^{-1}|D_d\Phi|^2
 &=u^{-1}\bigl(|\nabla\Phi|^2-|D_d\Phi|^2\bigr)-u^{-2}\ip{c(\nabla u)D_d\Phi}{\Phi}\\
 &\quad-\frac{d}{2}u^{-2}\nabla u\cdot\nabla|\Phi|^2+\frac{d(d-1)}{4}u^{-3}|\nabla u|^2|\Phi|^2.
\end{align*}
On the other hand, since $D_d^2=-\Delta$ and $-\ip{\Delta\Phi}{\Phi}=|\nabla\Phi|^2-\frac12\Delta|\Phi|^2$, we compute
\[
 \Div\bigl(u^{-1}\ip{c\,D_d\Phi}{\Phi}\bigr)
 =-u^{-2}\ip{c(\nabla u)D_d\Phi}{\Phi}+u^{-1}\Bigl(|\nabla\Phi|^2-|D_d\Phi|^2-\frac12\Delta|\Phi|^2\Bigr).
\]
Subtracting this divergence from the previous display, we are left with the scalar expression, written in terms of $f=u^{1-d}|\Phi|^2$,
\[
 \frac12u^{-1}\Delta(u^{d-1}f)-\frac{d}{2}u^{-2}\nabla u\cdot\nabla(u^{d-1}f)+\frac{d(d-1)}{4}u^{d-4}|\nabla u|^2f.
\]
Using the product rule and the identity
\[
 \zeta\Delta\zeta=\frac{d-2}{2}u^{d-3}\Bigl(\Delta u+\frac{d-4}{2}u^{-1}|\nabla u|^2\Bigr),
\]
we obtain $\frac12\Div(u^{d-2}\nabla f)+\frac{d-1}{d-2}f\,\zeta\Delta\zeta$. This proves \eqref{eq:sl}.
\end{proof}

Next we establish local $C^{2,\alpha}$ regularity and decay at infinity for solutions of \eqref{eq:coupled}. Set $Z=(U,V)$ and $\rho=|Z|=(|U|^2+|V|^2)^{1/2}$.

\begin{lemma}\label{lem:reg}
Let $d\ge3$. Every solution $Z=(U,V)\in\X_d^2$ of \eqref{eq:coupled} belongs to $C^{2,\alpha}_{\mathrm{loc}}(\R^d;\C^{2N_d})$ for all $0<\alpha<\nu$, and there exists a constant $C>0$ such that, for all $x\in\R^d$,
\begin{equation}\label{eq:decay}
 |Z(x)|\le C(1+|x|)^{1-d},\qquad |\nabla Z(x)|\le C(1+|x|)^{-d}.
\end{equation}
\end{lemma}

\begin{proof}
Set $\widetilde T(f,g)=(T_dg,T_df)$ and $M=\operatorname{diag}(|U|^\nu I_{N_d},|V|^\nu I_{N_d})$. Then $Z=\widetilde T(MZ)$ and $|M|\le\rho^\nu\in L^d(\R^d)$. Suppose $Z\in L^r$, starting with $r=q$, and choose a finite $s>r$ with $1/s+1/d\ge1/r$. Split $M=M_b+M_s$, with $M_b$ bounded and compactly supported and $\norm{M_s}_d$ small enough that $W\mapsto\widetilde T(M_sW)$ has norm less than one on both $L^r$ and $L^s$. These bounds follow from H\"older's inequality and fractional integration, with input exponents $1/(1/r+1/d)$ and $1/(1/s+1/d)$. Both exponents are greater than one and at most $r$, so $\widetilde T(M_bZ)$ belongs to both spaces. The Neumann series converges in both spaces to the same distribution. Uniqueness in $L^r$ identifies its sum with $Z$, hence $Z\in L^s$. Finitely many such steps give every finite exponent above $q$. In particular $MZ\in L^p\cap L^t$ for some $t>d$. The kernel of $\widetilde T$ is bounded by $C|x|^{1-d}$, which lies in $L^{t'}$ inside the unit ball and $L^{p'}$ outside it. Thus $Z\in L^\infty(\R^d)$.

By interpolation, $|U|^\nu U$ and $|V|^\nu V$ belong to $L^r(\R^d)$ for every $r\in[p,\infty)$. Representing $U=T_d(|V|^\nu V)$ and $V=T_d(|U|^\nu U)$, we deduce from Calder\'on--Zygmund bounds on the Riesz transforms in $\partial_jT_d$ that $\nabla Z\in L^r(\R^d)$ for all $r\in[p,\infty)$. Since $0<\nu\le1$, the map $z\mapsto|z|^\nu z$ lies in $C^1(\C^{N_d})$ with locally $\nu$-H\"older differential, so $|U|^\nu U,|V|^\nu V\in W^{1,r}_{\mathrm{loc}}(\R^d)$ for all $r<\infty$.

Because $-\Delta U=D_d(|V|^\nu V)$ and $-\Delta V=D_d(|U|^\nu U)$, interior $L^r$ elliptic estimates \cite{Taylor} imply $Z\in W^{2,r}_{\mathrm{loc}}(\R^d)$ for all $r<\infty$. By Morrey's embedding, $Z$ is locally Lipschitz; hence $|U|^\nu U,|V|^\nu V\in C^{1,\alpha}_{\mathrm{loc}}(\R^d)$ for every $\alpha\in(0,\nu)$, and interior Schauder estimates \cite{Taylor} show that $Z\in C^{2,\alpha}_{\mathrm{loc}}(\R^d)$.

Finally, by Lemma~\ref{lem:conf}, the pair $(U_{\mathcal K},-V_{\mathcal K})$ also solves \eqref{eq:coupled} in $\X_d^2$ and therefore belongs to $C^{2,\alpha}(B_1(0))$. Differentiating $W(x)=-|x|^{-d}c(x)W_{\mathcal K}(x/|x|^2)$ for $|x|\ge1$, we obtain \eqref{eq:decay} outside $B_1(0)$, and boundedness in $C^1(\overline{B_1(0)})$ covers $|x|\le1$.
\end{proof}

By Lemma~\ref{lem:reg}, the nonvanishing set
\[
 \Omega=\bigl\{x\in\R^d:\rho(x)>0\bigr\}
\]
is open in $\R^d$. On $\Omega$ we set
\[
 \begin{gathered}
  a=\frac{|U|^2}{\rho^2},\qquad \lambda=\frac{|U|^q+|V|^q}{\rho^2},\qquad d\mu=\lambda^d\,dx,\\
  w(a)=w_d(a)=\Bigl(a^{d/(d-1)}+(1-a)^{d/(d-1)}\Bigr)^{1-d},
 \end{gathered}
\]
extending $\lambda$ by $0$ on $\R^d\setminus\Omega$, and put $V_0=(d/2)^d|\Sph^d|$. Note that $1\le w(a)\le2$, $\rho^2=\lambda^{d-1}w(a)$, and $|U|^q+|V|^q=w(a)\lambda^d$, so
\begin{equation}\label{eq:I-vol}
 2I=\int_\Omega w(a)\,d\mu,\qquad I\le\mu(\Omega)\le2I.
\end{equation}
Since $q=2+\nu>2$, the function $z\mapsto|z|^q$ belongs to $C^{2,\nu}(\C^{N_d})$. On $\Omega$, where $\rho>0$ and $|U|^q+|V|^q>0$, Lemma~\ref{lem:reg} implies that both $\lambda$ and the composite function $w(a)=\rho^{2d}(|U|^q+|V|^q)^{1-d}$ belong to $C^{2,\alpha}(\Omega)$, and
\begin{equation}\label{eq:zeta-grad-bound}
 |\nabla a|+\lambda^{-1}|\nabla\lambda|\le C\rho^{-1}|\nabla Z|,\qquad \bigl|\nabla\lambda^{(d-2)/2}\bigr|^2\le C\rho^{-\nu}|\nabla Z|^2\qquad\text{on }\Omega.
\end{equation}

Following \cite[Proposition~7]{FL}, we now prove a global bound on $\rho^{-\nu}|\nabla Z|^2$ by applying Lemma~\ref{lem:sl} with $u=\rho_\eps^\nu$, where $\rho_\eps=(\rho^2+\eps^2)^{1/2}$.

\begin{lemma}
Let $d\ge3$. Every nontrivial solution $Z=(U,V)\in\X_d^2$ of \eqref{eq:coupled} satisfies
\begin{equation}\label{eq:weighted-grad}
 \int_\Omega\rho^{-\nu}|\nabla Z|^2\,dx\le C_d\int_{\R^d}\rho^q\,dx<\infty.
\end{equation}
Consequently, $\lambda^{(d-2)/2}$ belongs to $\dot H^1(\R^d)$ with weak gradient $\chi_\Omega\nabla\lambda^{(d-2)/2}$, where $\chi_\Omega$ is the indicator of $\Omega$, and
\begin{equation}\label{eq:conformal-sobolev}
 \frac{d-1}{d-2}\int_\Omega\bigl|\nabla\lambda^{(d-2)/2}\bigr|^2\,dx\ge\frac{d-1}{d}\,V_0^{2/d}\mu(\Omega)^{(d-2)/d}.
\end{equation}
\end{lemma}

\begin{proof}
Fix $\eps>0$ and set $\rho_\eps=(\rho^2+\eps^2)^{1/2}\in C^{2,\alpha}_{\mathrm{loc}}(\R^d)$ and $\zeta_\eps=\rho_\eps^{(d-2)/(d-1)}$. We apply Lemma~\ref{lem:sl} on $\R^d$ with $u=\rho_\eps^\nu$ to $\Phi=U$ and $\Phi=V$ and sum the two equations:
\begin{align*}
 \rho_\eps^2\Bigl(\bigl|\T(\rho_\eps^{-d/(d-1)}U)\bigr|^2+\bigl|\T(\rho_\eps^{-d/(d-1)}V)\bigr|^2\Bigr)-\frac{d-1}{d-2}\frac{\rho^2}{\rho_\eps^2}\zeta_\eps\Delta\zeta_\eps
 &=\frac{d-1}{d}\rho_\eps^{-\nu}\bigl(|D_dU|^2+|D_dV|^2\bigr)\\
 &\quad+\Div\mathcal V_\eps,
\end{align*}
where $\mathcal V_\eps=\rho_\eps^{-\nu}(\ip{c\,D_dU}{U}+\ip{c\,D_dV}{V})+\frac12\rho_\eps^{\nu(d-2)}\nabla(\rho^2/\rho_\eps^2)$. By \eqref{eq:decay}, both $|\mathcal V_\eps(x)|$ and $|\nabla\zeta_\eps(x)|$ are $O(|x|^{1-2d})$ as $|x|\to\infty$, so integrating over $B_R(0)$ and letting $R\to\infty$ eliminates the boundary terms. Integrating by parts in the $\zeta_\eps\Delta\zeta_\eps$ term, we find
\begin{align*}
 -\int_{\R^d}\frac{\rho^2}{\rho_\eps^2}\zeta_\eps\Delta\zeta_\eps\,dx
 &=\int_{\R^d}\frac{\rho^2}{\rho_\eps^2}|\nabla\zeta_\eps|^2\,dx+\int_{\R^d}\zeta_\eps\nabla\Bigl(1-\frac{\eps^2}{\rho_\eps^2}\Bigr)\cdot\nabla\zeta_\eps\,dx\\
 &\ge\Bigl(\frac{d-2}{d-1}\Bigr)^2\int_{\R^d}\frac{\rho^2}{\rho_\eps^2}\rho_\eps^{-\nu}|\nabla\rho_\eps|^2\,dx,
\end{align*}
since $\nabla(1-\eps^2/\rho_\eps^2)$ and $\nabla\zeta_\eps$ are nonnegative multiples of $\nabla\rho_\eps$. Using $|D_dU|^2+|D_dV|^2\le\rho^{2+2\nu}\le\rho_\eps^\nu\rho^q$, we obtain
\begin{equation}\label{eq:eps-apriori}
 \begin{aligned}
  &\int_{\R^d}\rho_\eps^2\Bigl(\bigl|\T(\rho_\eps^{-d/(d-1)}U)\bigr|^2+\bigl|\T(\rho_\eps^{-d/(d-1)}V)\bigr|^2\Bigr)\,dx\\
  &\qquad+\frac{d-2}{d-1}\int_{\R^d}\frac{\rho^2}{\rho_\eps^2}\rho_\eps^{-\nu}|\nabla\rho_\eps|^2\,dx
  \le\frac{d-1}{d}\int_{\R^d}\rho^q\,dx.
 \end{aligned}
\end{equation}
With $\varphi_U=\rho_\eps^{-d/(d-1)}U$, the product rule and \eqref{eq:twistor-norm} imply $\rho_\eps^2|\nabla\varphi_U|^2=\rho_\eps^2|\T\varphi_U|^2+\frac1d\rho_\eps^2|D_d\varphi_U|^2$ together with
\begin{gather*}
 \rho_\eps^{-\nu}|\nabla U|^2\le2\rho_\eps^2|\nabla\varphi_U|^2+\frac{2d^2}{(d-1)^2}\frac{|U|^2}{\rho_\eps^2}\rho_\eps^{-\nu}|\nabla\rho_\eps|^2,\\
 \rho_\eps^2|D_d\varphi_U|^2\le2\rho_\eps^{-\nu}|D_dU|^2+\frac{2d^2}{(d-1)^2}\frac{|U|^2}{\rho_\eps^2}\rho_\eps^{-\nu}|\nabla\rho_\eps|^2.
\end{gather*}
We add the analogous bounds for $V$, integrate over $\R^d$ using \eqref{eq:eps-apriori}, and let $\eps\downarrow0$ by monotone convergence to obtain \eqref{eq:weighted-grad}.

By \eqref{eq:zeta-grad-bound} and \eqref{eq:weighted-grad}, $\chi_\Omega\nabla\lambda^{(d-2)/2}\in L^2(\R^d)$. For each $\delta>0$, the truncation $(\lambda^{(d-2)/2}-\delta)_+$ is Lipschitz and compactly supported in $\Omega$ by \eqref{eq:decay}; letting $\delta\downarrow0$ shows that $\lambda^{(d-2)/2}\in\dot H^1(\R^d)$. By the sharp Sobolev inequality on $\R^d$ \cite{Aubin,Talenti,Beckner,FL-HLS},
\[
 \int_\Omega\bigl|\nabla\lambda^{(d-2)/2}\bigr|^2\,dx\ge\frac{d(d-2)}{4}|\Sph^d|^{2/d}\mu(\Omega)^{(d-2)/d},
\]
and \eqref{eq:conformal-sobolev} follows since $|\Sph^d|^{2/d}=(2/d)^2V_0^{2/d}$.
\end{proof}

For the components $Y$ of a $C^1$ vector field in the orthonormal frame $(\lambda^{-1}\partial_j)_{j=1}^d$ on $\Omega$, write
\[
 \Divmu Y=\lambda^{-d}\Div\bigl(\lambda^{d-1}Y\bigr),\qquad \norm{Y}_{L^2(\Omega,d\mu)}^2=\int_\Omega|Y|^2\,d\mu.
\]

\begin{lemma}\label{lem:ibp}
Let $d\ge3$. Every vector field $Y\in C^1(\Omega;\R^d)\cap L^2(\Omega,d\mu)$ with $\Divmu Y\in L^1(\Omega,d\mu)$ satisfies
\begin{equation}\label{eq:div-zero}
 \int_\Omega\Divmu Y\,d\mu=0.
\end{equation}
\end{lemma}

\begin{proof}
Fix $\chi\in C^\infty(\R;[0,1])$ equal to $1$ on $(-\infty,1]$ and $0$ on $[2,\infty)$, and put $\theta_\eps(x)=1-\chi(\rho(x)/\eps)$ and $\chi_R(x)=\chi(|x|/R)$ for $\eps>0$ and $R>1$. Since $\theta_\eps\chi_R\in C_c^1(\Omega)$, the divergence theorem reads
\[
 \int_\Omega\theta_\eps\chi_R\,\Divmu Y\,d\mu=-\int_\Omega\chi_R(\lambda^{-1}\nabla\theta_\eps)\cdot Y\,d\mu-\int_\Omega\theta_\eps(\lambda^{-1}\nabla\chi_R)\cdot Y\,d\mu.
\]
On $\{\eps<\rho<2\eps\}$, we have $|\nabla\theta_\eps|\le C\rho^{-1}|\nabla Z|$ and $\lambda^{d-2}\rho^{-2}\le\rho^{-\nu}$, so by \eqref{eq:weighted-grad},
\[
 \int_\Omega|\lambda^{-1}\nabla\theta_\eps|^2\,d\mu\le C\int_{\{\eps<\rho<2\eps\}}\rho^{-\nu}|\nabla Z|^2\,dx\longrightarrow0\qquad\text{as }\eps\downarrow0.
\]
On $A_R=\{R<|x|<2R\}\cap\Omega$, Lemma~\ref{lem:reg} implies $\lambda\le\rho^\nu\le CR^{-2}$, and hence
\[
 \int_{A_R}|\lambda^{-1}\nabla\chi_R|^2\,d\mu\le CR^{-2}\int_{A_R}\lambda^{d-2}\,dx\le CR^{2-d}\longrightarrow0\qquad\text{as }R\to\infty.
\]
By Cauchy--Schwarz in $L^2(\Omega,d\mu)$, letting $\eps\downarrow0$ and then $R\to\infty$ proves \eqref{eq:div-zero}.
\end{proof}

\subsection{Volume and mass identities}
On $\Omega$, we introduce the rescaled spinors $(H,J)$ and the auxiliary $C^{2,\alpha}(\Omega)$ spinors $(\Phi_1,\Phi_2)$ by
\[
 \begin{gathered}
  H=\lambda^{-1/2}\frac{U}{\rho},\qquad J=\lambda^{-1/2}\frac{V}{\rho},\\
  \Phi_1=w(a)^{-1/2}U=\lambda^{d/2}H,\qquad \Phi_2=w(a)^{-1/2}V=\lambda^{d/2}J.
 \end{gathered}
\]

Define $B,X,T,\Theta$ by
\[
 \begin{gathered}
  B_j=\frac{\ip{c_jV}{U}}{\rho^2}=\lambda\ip{c_jJ}{H},\qquad X=\lambda^{-1}\nabla a+\frac2dB,\\
  T=\lambda^{-1}\bigl(|\T H|^2+|\T J|^2\bigr),\qquad \Theta=|X|^2+\frac{4}{d^2}\bigl(a(1-a)-|B|^2\bigr).
 \end{gathered}
\]
Set
\begin{equation}\label{eq:gamma-def}
 \gamma(a)=\gamma_d(a)=\frac{w'(a)}{2d\,w(a)}=\frac{(1-a)^{1/(d-1)}-a^{1/(d-1)}}{2\bigl(a^{d/(d-1)}+(1-a)^{d/(d-1)}\bigr)}\qquad(a\in[0,1]).
\end{equation}
Note that $|H|^2=a\lambda^{-1}$, $|J|^2=(1-a)\lambda^{-1}$, and $2d\gamma(a)\nabla a=w(a)^{-1}\nabla(w(a))\in C^{1,\alpha}(\Omega)$. Wherever $J\ne0$, the spinors $J,c_1J,\dots,c_dJ$ are pairwise orthogonal of squared norm $(1-a)\lambda^{-1}$, so Bessel's inequality implies
\begin{equation}\label{eq:bessel-B}
 |B|^2+\lambda^2\ip{H}{J}^2\le a(1-a)\qquad\text{on }\Omega.
\end{equation}

\begin{lemma}
Let $d\ge3$. Everywhere on $\Omega$, the spinors $\Phi_1,\Phi_2\in C^{2,\alpha}(\Omega)$ solve
\begin{equation}\label{eq:D-Phi}
 \begin{aligned}
  D_d\Phi_1&=-d\gamma(a)\,c(\nabla a)\Phi_1+\lambda\bigl(1+2a\gamma(a)\bigr)\Phi_2,\\
  D_d\Phi_2&=-d\gamma(a)\,c(\nabla a)\Phi_2+\lambda\bigl(1-2(1-a)\gamma(a)\bigr)\Phi_1,
 \end{aligned}
\end{equation}
we have
\begin{equation}\label{eq:TH-H}
 \ip{\T_jH}{H}=\frac{1+2a\gamma(a)}{2}\,X_j,\qquad \ip{\T_jJ}{J}=-\frac{1-2(1-a)\gamma(a)}{2}\,X_j,
\end{equation}
and the vector fields $B$ and $X$ satisfy
\begin{equation}\label{eq:BX-flux}
 \begin{gathered}
  \Divmu\bigl(w(a)B\bigr)=-\bigl(1-2a+4a(1-a)\gamma(a)\bigr)w(a),\\
  \Divmu\bigl(w(a)X\bigr)=2w(a)\lambda^{-1}\bigl((1-a)|\T H|^2-a|\T J|^2\bigr).
 \end{gathered}
\end{equation}
\end{lemma}

\begin{proof}
By \eqref{eq:gamma-def}, we have
\[
 \frac{(1-a)^{1/(d-1)}}{a^{d/(d-1)}+(1-a)^{d/(d-1)}}=1+2a\gamma(a),\qquad
 \frac{a^{1/(d-1)}}{a^{d/(d-1)}+(1-a)^{d/(d-1)}}=1-2(1-a)\gamma(a),
\]
so $w(a)^{-1/2}|V|^{q-2}V=\lambda(1+2a\gamma(a))\Phi_2$ and $w(a)^{-1/2}|U|^{q-2}U=\lambda(1-2(1-a)\gamma(a))\Phi_1$ on $\Omega$. Differentiating $\Phi_1=w(a)^{-1/2}U$ and $\Phi_2=w(a)^{-1/2}V$ via $\nabla(w(a))=2d\gamma(a)w(a)\nabla a$ and \eqref{eq:coupled}, we obtain \eqref{eq:D-Phi}. Computing $\T_j(\lambda^{-d/2}\Phi_1)$ and $\T_j(\lambda^{-d/2}\Phi_2)$ as in Lemma~\ref{lem:sl}, we see that the $\nabla\lambda$ terms cancel in $\ip{\T_jH}{H}$ and $\ip{\T_jJ}{J}$, leaving \eqref{eq:TH-H}.

Next, since $\lambda^{d-1}w(a)B_j=\ip{c_jV}{U}$, it follows from \eqref{eq:coupled} that
\[
 \Div\bigl(\lambda^{d-1}w(a)B\bigr)=\ip{D_dV}{U}-\ip{V}{D_dU}=|U|^q-|V|^q=-\bigl(1-2a+4a(1-a)\gamma(a)\bigr)w(a)\lambda^d,
\]
which is the first identity in \eqref{eq:BX-flux}. From \eqref{eq:D-Phi}, we compute
\begin{align*}
 \lambda^{-1}\ip{c\,D_d\Phi_1}{\Phi_1}&=\lambda^{d-1}\bigl(d\gamma(a)a\,\lambda^{-1}\nabla a+(1+2a\gamma(a))B\bigr),\\
 \lambda^{-1}\ip{c\,D_d\Phi_2}{\Phi_2}&=\lambda^{d-1}\bigl(d\gamma(a)(1-a)\lambda^{-1}\nabla a-(1-2(1-a)\gamma(a))B\bigr),
\end{align*}
and hence
\begin{equation}\label{eq:flux-combos}
 \begin{gathered}
  \lambda^{-1}\bigl(\ip{c\,D_d\Phi_1}{\Phi_1}+\ip{c\,D_d\Phi_2}{\Phi_2}\bigr)=d\gamma(a)\lambda^{d-1}X,\\
  (1-a)\lambda^{-1}\ip{c\,D_d\Phi_1}{\Phi_1}-a\lambda^{-1}\ip{c\,D_d\Phi_2}{\Phi_2}=\lambda^{d-1}B,
 \end{gathered}
\end{equation}
together with
\begin{equation}\label{eq:DPhi-norms}
 \begin{gathered}
  \lambda^{-1}\bigl(|D_d\Phi_1|^2+|D_d\Phi_2|^2\bigr)=\lambda^d\bigl(1+d^2\gamma(a)^2\Theta\bigr),\\
  (1-a)\lambda^{-1}|D_d\Phi_1|^2-a\lambda^{-1}|D_d\Phi_2|^2=-\lambda^d\Divmu B.
 \end{gathered}
\end{equation}
We now apply Lemma~\ref{lem:sl} on $\Omega$ with $u=\lambda$ to $\Phi_1$ and $\Phi_2$. Taking $(1-a)$ times the first equation minus $a$ times the second cancels $\zeta\Delta\zeta$. The weighted difference of the full fluxes is $\lambda^{d-1}(X/2+(d-1)B/d)$, while their sum is $d\gamma\lambda^{d-1}X$. Moving the weights inside the divergences therefore adds $d\gamma\lambda^{d-1}\nabla a\cdot X$. By \eqref{eq:DPhi-norms}, the $\Divmu B$ terms cancel, leaving
\[
 \Divmu X+2d\gamma(a)(\lambda^{-1}\nabla a)\cdot X=2\lambda^{-1}\bigl((1-a)|\T H|^2-a|\T J|^2\bigr).
\]
Multiplying by $w(a)$, we obtain the second identity in \eqref{eq:BX-flux}.
\end{proof}

Set on $\Omega$
\[
 \mathcal E_0=T-d(d-1)\gamma(a)^2\Theta.
\]
By \eqref{eq:zeta-grad-bound} and the definitions of $H,J$, we have
\[
 \lambda^d\bigl(T+|\lambda^{-1}\nabla a|^2\bigr)
 \le C\lambda^{d-2}\rho^{-2}|\nabla Z|^2
 \le C\rho^{-\nu}|\nabla Z|^2.
\]
Thus \eqref{eq:weighted-grad}, the boundedness of $B,\gamma$, and $\mu(\Omega)<\infty$ give $X\in L^2(d\mu)$ and $T,\Theta,\mathcal E_0\in L^1(d\mu)$.

\begin{proposition}
Let $d\ge3$. Every nontrivial solution $(U,V)\in\X_d^2$ of \eqref{eq:coupled} satisfies
\begin{equation}\label{eq:vol-id}
 \frac{d-1}{d}\mu(\Omega)=\frac{d-1}{d-2}\int_\Omega\bigl|\nabla\lambda^{(d-2)/2}\bigr|^2\,dx+\int_\Omega\mathcal E_0\,d\mu.
\end{equation}
\end{proposition}

\begin{proof}
Set $\zeta=\lambda^{(d-2)/2}$ on $\Omega$. We add the two instances of Lemma~\ref{lem:sl} at $u=\lambda$ for $\Phi_1$ and $\Phi_2$. Substituting $\zeta\Delta\zeta=\Div(\zeta\nabla\zeta)-|\nabla\zeta|^2$ and using \eqref{eq:flux-combos} and \eqref{eq:DPhi-norms}, we obtain on $\Omega$
\begin{equation}\label{eq:vol-pointwise}
 \frac{d-1}{d-2}|\nabla\zeta|^2+\lambda^d\mathcal E_0-\frac{d-1}{d}\lambda^d=\lambda^d\Divmu Y,\qquad Y=\frac{d-1}{2}\lambda^{-2}\nabla\lambda+d\gamma(a)X.
\end{equation}
By \eqref{eq:flux-combos}, $\lambda^{d-1}Y=\frac{d-1}{d-2}\zeta\nabla\zeta+\lambda^{-1}(\ip{c\,D_d\Phi_1}{\Phi_1}+\ip{c\,D_d\Phi_2}{\Phi_2})$ belongs to $C^{1,\alpha}(\Omega;\R^d)$. Since $\nabla\zeta\in L^2(dx)$ and $X\in L^2(d\mu)$, we have $Y\in L^2(d\mu)$, and \eqref{eq:vol-pointwise} gives $\Divmu Y\in L^1(d\mu)$. Integrating \eqref{eq:vol-pointwise} over $\Omega$ via Lemma~\ref{lem:ibp} proves \eqref{eq:vol-id}.
\end{proof}

To relate $\mu(\Omega)$ to $2I=\int_\Omega w(a)\,d\mu$, use
\[
 (\lambda^{-1}\nabla a)\cdot\Bigl(B+\frac{d}{2}X\Bigr)=\frac{d}{2}\Theta-\frac{2}{d}a(1-a).
\]
For any $h\in C^1[0,1]$, writing $h(a)(B+\frac{d}{2}X)=(h(a)/w(a))\,w(a)(B+\frac{d}{2}X)$ and using \eqref{eq:BX-flux}, we compute on $\Omega$
\begin{equation}\label{eq:h-div-raw}
 \begin{aligned}
  \Divmu\biggl(h(a)\Bigl(B+\frac{d}{2}X\Bigr)\biggr)
  &=\Bigl(-(1-2a)h(a)-\frac{2}{d}a(1-a)h'(a)\Bigr)\\
  &\quad+d\,h(a)\lambda^{-1}\bigl((1-a)|\T H|^2-a|\T J|^2\bigr)+\frac{d}{2}\bigl(h'-2d\gamma h\bigr)(a)\Theta.
 \end{aligned}
\end{equation}
We choose $h$ so that the first line on the right is $w(a)$ minus a constant. The requirement that $h$ be bounded at both endpoints determines the constant.

\begin{lemma}\label{lem:h-ode}
For every $d\ge3$, the ratio $\omega_d=2S_d^d/V_0$ satisfies $1<\omega_d<2$ and
\begin{equation}\label{eq:omega-mean}
 \omega_d=\frac{\int_0^1[s(1-s)]^{d/2-1}w(s)\,ds}{\int_0^1[s(1-s)]^{d/2-1}\,ds},
\end{equation}
with $\omega_3=16/(3\pi)$ and $\omega_4=9\pi/16$. Consequently, the linear equation
\begin{equation}\label{eq:h-ode}
 -(1-2a)h_d(a)-\frac{2}{d}a(1-a)h_d'(a)=w(a)-\omega_d
\end{equation}
admits a unique bounded solution $h_d$ on $(0,1)$, given by
\begin{equation}\label{eq:h-int}
 h_d(a)=\frac{d}{2[a(1-a)]^{d/2}}\int_0^a[s(1-s)]^{d/2-1}\bigl(\omega_d-w(s)\bigr)\,ds,
\end{equation}
which extends to $h_d\in C^{1,1/(d-1)}[0,1]$ and satisfies $h_d(1-a)=-h_d(a)$ and $h_d(0)=\omega_d-1$. Moreover, if $w$ is concave on $[0,1]$, then $h_d'<0$ on $(0,1/2]$ and $0<h_d<\omega_d-1$ on $(0,1/2)$.
\end{lemma}

\begin{proof}
On the pair $(U_\eta,V_\eta)$ from Lemma~\ref{lem:trial}, we have $\lambda(x)=d/(1+|x|^2)$, $a(x)=(1+|x|^2)^{-1}$, $\mu(\R^d)=V_0$, and $I=S_d^d$. Passing to polar coordinates and substituting $s=(1+r^2)^{-1}$ in $\mu(\R^d)$ and in $2I=\int_{\R^d}w(a)\,d\mu$, we obtain \eqref{eq:omega-mean}. Since $1\le w\le2$ and $w$ is nonconstant, $1<\omega_d<2$. The values of $\omega_3$ and $\omega_4$ follow from $S_3^3=18\pi$ and $S_4^4=12\pi^3$.

Multiplying \eqref{eq:h-ode} by $[a(1-a)]^{d/2-1}$ turns the left-hand side into $-\frac{2}{d}\frac{\mathrm d}{\mathrm da}\bigl([a(1-a)]^{d/2}h_d(a)\bigr)$. By \eqref{eq:omega-mean}, the integral in \eqref{eq:h-int} vanishes at $a=1$, so \eqref{eq:h-int} is the unique solution bounded near $0$ and $1$. Symmetry $w(1-s)=w(s)$ implies $h_d(1-a)=-h_d(a)$. Substituting $s=a\tau$ in \eqref{eq:h-int} on $(0,1/2]$, we write
\[
 h_d(a)=\frac{d}{2(1-a)^{d/2}}\int_0^1\tau^{d/2-1}(1-a\tau)^{d/2-1}\bigl(\omega_d-w(a\tau)\bigr)\,d\tau.
\]
Since $w\in C^{1,1/(d-1)}[0,1/2]$ with $w(0)=1$, differentiating under the integral sign shows that $h_d\in C^{1,1/(d-1)}[0,1]$ and $h_d(0)=\omega_d-1$.

Finally, since $[a(1-a)]^{d/2}=\frac{d}{2}\int_0^a[s(1-s)]^{d/2-1}(1-2s)\,ds$, combining \eqref{eq:h-ode} and \eqref{eq:h-int} we obtain, for $a\in(0,1)$,
\begin{equation}\label{eq:hd-prime-int}
 \begin{aligned}
  -\frac{4}{d^2}[a(1-a)]^{d/2+1}h_d'(a)
  &=\int_0^a[s(1-s)]^{d/2-1}\Bigl[2(a-s)\bigl(w(a)-\omega_d\bigr)\\
  &\qquad\qquad\qquad\qquad\quad+(1-2a)\bigl(w(a)-w(s)\bigr)\Bigr]\,ds.
 \end{aligned}
\end{equation}
If $w$ is concave on $[0,1]$, then $w(a)-w(s)\ge2d\gamma(a)w(a)(a-s)$ for $0\le s\le a\le1/2$, while concavity at $a\le1/2$ also implies $w(a)[1+d(1-2a)\gamma(a)]=w(a)+w'(a)(\frac12-a)\ge w(1/2)=2>\omega_d$. By \eqref{eq:hd-prime-int}, $h_d'<0$ on $(0,1/2]$, and since $h_d(1/2)=0$ we have $0<h_d<\omega_d-1$ on $(0,1/2)$.
\end{proof}

With $h_d$ defined by Lemma~\ref{lem:h-ode}, set on $\Omega$
\begin{equation}\label{eq:E1-def}
 \begin{aligned}
  \mathcal E_1&=\lambda^{-1}\Bigl(1-\frac{d-1}{\omega_d}(1-a)h_d(a)\Bigr)|\T H|^2+\lambda^{-1}\Bigl(1+\frac{d-1}{\omega_d}a\,h_d(a)\Bigr)|\T J|^2\\
  &\quad-\Bigl[\frac{d-1}{2\omega_d}\bigl(h_d'-2d\gamma h_d\bigr)(a)+d(d-1)\gamma(a)^2\Bigr]\Theta.
 \end{aligned}
\end{equation}

\begin{proposition}
Let $d\ge3$. Every nontrivial solution $(U,V)\in\X_d^2$ of \eqref{eq:coupled} satisfies
\begin{equation}\label{eq:second-id}
 \frac{d-1}{d}\Bigl(\frac{I}{S_d^d}V_0-\mu(\Omega)\Bigr)=\int_\Omega(\mathcal E_1-\mathcal E_0)\,d\mu.
\end{equation}
\end{proposition}

\begin{proof}
Since $X\in L^2(d\mu)$, $B,h_d$ are bounded and $\mu(\Omega)<\infty$, the field $Y=h_d(a)(B+\frac{d}{2}X)$ belongs to $C^1(\Omega;\R^d)\cap L^2(d\mu)$. As $T,\Theta\in L^1(d\mu)$, \eqref{eq:h-div-raw} gives $\Divmu Y\in L^1(d\mu)$. Integrating by Lemma~\ref{lem:ibp} and using \eqref{eq:h-ode} gives
\[
 \int_\Omega(\mathcal E_1-\mathcal E_0)\,d\mu=\frac{d-1}{d\omega_d}\int_\Omega\bigl(w(a)-\omega_d\bigr)\,d\mu.
\]
Substituting $\int_\Omega w(a)\,d\mu=2I$ and $\omega_dV_0=2S_d^d$ gives \eqref{eq:second-id}.
\end{proof}

\section{Mass bounds}

We combine the volume and mass identities with the sharp Sobolev inequality, and then bound the twistor terms in terms of $|X|^2$.

\begin{proposition}\label{prop:kappa-red}
Let $d\ge3$ and let $(U,V)\in\X_d^2$ be a nontrivial solution of \eqref{eq:coupled} with $\int_\Omega\mathcal E_0\,d\mu\ge0$. Then $s=(\mu(\Omega)/V_0)^{1/d}\ge1$ and for every $\kappa\ge0$ we have
\begin{equation}\label{eq:kappa-red}
 \frac{d-1}{d}\Bigl(\frac{I}{S_d^d}-1\Bigr)V_0\ge\frac{d-1}{d}\Bigl((1+\kappa)s^{d-2}-\kappa s^d-1\Bigr)V_0+\int_\Omega\bigl(\mathcal E_1+\kappa\mathcal E_0\bigr)\,d\mu.
\end{equation}
Moreover, if $I\le S_d^d$, then $1\le s^d\le\omega_d<2$.
\end{proposition}

\begin{proof}
By \eqref{eq:vol-id} and \eqref{eq:conformal-sobolev},
\[
 0\le\int_\Omega\mathcal E_0\,d\mu\le\frac{d-1}{d}V_0(s^d-s^{d-2}).
\]
Since $s>0$ by \eqref{eq:I-vol}, we have $s\ge1$. In \eqref{eq:second-id}, write $\mathcal E_1-\mathcal E_0=(\mathcal E_1+\kappa\mathcal E_0)-(1+\kappa)\mathcal E_0$ and apply the upper bound just obtained. This gives \eqref{eq:kappa-red}. If $I\le S_d^d$, then $\mu(\Omega)\le2I\le2S_d^d=\omega_dV_0$ by \eqref{eq:I-vol}, so $s^d\le\omega_d<2$.
\end{proof}

\begin{lemma}\label{lem:twistor-lower}
Let $d\ge3$. On $\{0<a<1\}\subset\Omega$, for every $\tau\in\R$ satisfying $1-(1-a)\tau\ge0$ and $1+a\tau\ge0$,
\begin{equation}\label{eq:twistor-weighted}
 \begin{aligned}
  &\lambda^{-1}\Bigl[\bigl(1-(1-a)\tau\bigr)|\T H|^2+(1+a\tau)|\T J|^2\Bigr]\\
  &\qquad\ge\frac{d}{4(d-1)}\Bigl(\frac{1-(1-2a)\tau}{a(1-a)}-4\tau\gamma(a)+4\gamma(a)^2\Bigr)|X|^2.
 \end{aligned}
\end{equation}
In particular, at $\tau=0$ we have $T\ge\frac{d}{4(d-1)}\bigl(\frac{1}{a(1-a)}+4\gamma(a)^2\bigr)|X|^2$.
\end{lemma}

\begin{proof}
For $\varphi\in\C^{N_d}$ and $1\le j\le d$, define $\Pi^{(j)}(\varphi)\in(\C^{N_d})^d$ by $\Pi^{(j)}(\varphi)_k=\delta_{jk}\varphi+\frac1dc_kc_j\varphi$, so that $\ip{\Pi^{(j)}(\varphi)}{\Pi^{(\ell)}(\varphi)}=\frac{d-1}{d}\delta_{j\ell}|\varphi|^2$. Since $\sum_{k=1}^dc_k\T_kH=0$, we have $\ip{\T H}{\Pi^{(j)}(H)}=\ip{\T_jH}{H}$. By Bessel's inequality along $\{\Pi^{(j)}(H)\}_{j=1}^d$ and along $\{\Pi^{(j)}(J)\}_{j=1}^d$ together with \eqref{eq:TH-H},
\[
 \lambda^{-1}|\T H|^2\ge\frac{d\bigl(1+2a\gamma(a)\bigr)^2}{4(d-1)a}|X|^2,\qquad
 \lambda^{-1}|\T J|^2\ge\frac{d\bigl(1-2(1-a)\gamma(a)\bigr)^2}{4(d-1)(1-a)}|X|^2.
\]
Multiplying the two inequalities by $1-(1-a)\tau$ and $1+a\tau$, respectively, and adding gives \eqref{eq:twistor-weighted}.
\end{proof}

\subsection{Dimension three}
Throughout this subsection $d=3$, $q=3$, $S_3^3=18\pi$, $V_0=27\pi^2/4$, and $\omega_3=16/(3\pi)$.

\begin{proposition}\label{prop:d3}
Every nontrivial solution $(U,V)\in\X_3^2$ of \eqref{eq:coupled} satisfies $I\ge18\pi$. If $I=18\pi$, then $T=0$ and $\Theta=0$ everywhere on $\Omega$.
\end{proposition}

\begin{proof}
Since $\dim_\R\C^2=4$, wherever $J\ne0$ on $\Omega$ the spinors $J,c_1J,c_2J,c_3J$ form a real orthogonal basis of $\C^2$ of squared norm $(1-a)\lambda^{-1}$ (cf.\ \cite[Lemma~8.2]{WZcurl}). With $t=\ip{U}{V}/\rho^2=\lambda\ip{H}{J}\in C^{2,\alpha}(\Omega)$, Parseval's identity reads
\[
 t^2+|B|^2=a(1-a)\qquad\text{on }\Omega.
\]
To control the negative term in $\mathcal E_0=T-6\gamma_3^2\bigl(|X|^2+\frac49t^2\bigr)$, we define on $\Omega$ the vector field $Y=\lambda^{-1}t\nabla t-4\gamma_3t^2B$. We claim that $Y\in C^{1,\alpha}(\Omega;\R^3)$ and
\begin{equation}\label{eq:d3-polarized}
 \begin{aligned}
  3\mathcal E_0+2\,\Divmu Y
  &=2\bigl|\lambda^{-1}\nabla t-\gamma_3t(3X+2B)\bigr|^2+\frac23\gamma_3^2t^2|3X+2B|^2+32\gamma_3^2t^4\\
  &\quad+(3-4t^2)T+4t\lambda^{-1}\ip{\T H}{\T J}-6\gamma_3^2(3+4t^2)|X|^2\\
   &\quad+\frac83|\gamma_3|t^2\Bigl(2|1-2a|-\bigl(3+4a(1-a)\bigr)|\gamma_3|\Bigr).
 \end{aligned}
\end{equation}
Apply Lemma~\ref{lem:sl} on $\Omega$ at $u=\lambda$ and $\zeta=\lambda^{1/2}$ to $\Phi_1+\Phi_2$ and $\Phi_1-\Phi_2$. We take one quarter of the difference and use \eqref{eq:D-Phi} to compute $\ip{D_3\Phi_1}{D_3\Phi_2}$. This gives
\[
 \lambda^2\ip{\T H}{\T J}-\frac23\Bigl[\bigl(1+2a\gamma_3\bigr)\bigl(1-2(1-a)\gamma_3\bigr)+9\gamma_3^2|\lambda^{-1}\nabla a|^2\Bigr]\lambda^3t=2t\zeta\Delta\zeta+\Div\mathcal J,
\]
where $\mathcal J\in C^{1,\alpha}(\Omega;\R^3)$ is given via \eqref{eq:D-Phi} by
\[
 \mathcal J=\frac12\lambda^{-1}\bigl(\ip{c\,D_3\Phi_1}{\Phi_2}+\ip{c\,D_3\Phi_2}{\Phi_1}\bigr)+\frac12\lambda\nabla t=3\gamma_3\lambda t\nabla a+\frac12\lambda\nabla t.
\]
To eliminate the curvature term, use \eqref{eq:vol-pointwise} in the form $2\zeta\Delta\zeta=\lambda^3(\mathcal E_0-\frac23)-\Div(3\gamma_3\lambda^2X)$. Since $3\gamma_3\lambda^2X\in C^{1,\alpha}(\Omega;\R^3)$ by \eqref{eq:flux-combos}, the field $Y=\lambda^{-2}(2t\mathcal J-6\gamma_3\lambda^2t^2X)$ belongs to $C^{1,\alpha}(\Omega;\R^3)$. Denote the coefficient in square brackets in the polarized equation by $K$. By the product rule,
\[
 \Divmu Y=2t\lambda^{-1}\ip{\T H}{\T J}+\frac43t^2(1-K)-2t^2\mathcal E_0+|\lambda^{-1}\nabla t|^2-2\gamma_3t(\lambda^{-1}\nabla t)\cdot(3X+2B).
\]
Using $\lambda^{-1}\nabla a=X-\frac23B$ and $|B|^2=a(1-a)-t^2$ gives
\[
 K-1=-2(1-2a)\gamma_3-4\gamma_3^2t^2+9\gamma_3^2|X|^2-12\gamma_3^2X\cdot B.
\]
Inserting this and the definition of $\mathcal E_0$, and completing the square in $\lambda^{-1}\nabla t$, gives \eqref{eq:d3-polarized}. In its last term we used $\gamma_3(1-2a)=|\gamma_3|\,|1-2a|$.

We claim that the right-hand side of \eqref{eq:d3-polarized} is nonnegative on $\Omega$. On $\{a(1-a)=0\}$, we have $t=0$ and $X=0$, so the right-hand side equals $3T+2|\lambda^{-1}\nabla t|^2\ge0$. On $\{0<a<1\}$, put $m=\sqrt{a(1-a)}\in(0,1/2]$. Then $|t|\le m$, $2\lambda^{-1}|\ip{\T H}{\T J}|\le T$, $|\gamma_3|\le\frac12\sqrt{1-4m^2}$ by Lemma~\ref{lem:d3-profiles}\textup{(i)} below, and $3|X|^2\le8m^2T$ by Lemma~\ref{lem:twistor-lower} at $\tau=0$. Hence the second line of \eqref{eq:d3-polarized} is bounded from below by $T$ times
\[
 3-2m-4m^2-4m^2(1-4m^2)(3+4m^2)=(1-2m)+2(1-4m^2)^2+64m^6>0.
\]
By \eqref{eq:d3-w-gamma} below, the last term on the right-hand side of \eqref{eq:d3-polarized} simplifies to $\frac83\gamma_3^2t^2(1-2m)(1+6m)\ge0$. Therefore $3\mathcal E_0+2\,\Divmu Y\ge0$ on $\Omega$. By \eqref{eq:zeta-grad-bound}, \eqref{eq:weighted-grad}, and \eqref{eq:d3-polarized}, we have $Y\in L^2(\Omega,d\mu)$ and $\Divmu Y\in L^1(\Omega,d\mu)$. Lemma~\ref{lem:ibp} therefore implies $\int_\Omega\Divmu Y\,d\mu=0$, and hence $\int_\Omega\mathcal E_0\,d\mu\ge0$.

By Lemma~\ref{lem:d3-profiles}\textup{(ii)} below, all three coefficients of $\mathcal E_1$ in \eqref{eq:E1-def} are strictly positive on $(0,1)$, and since $\Theta=0$ on $\{a(1-a)=0\}$ we have $\mathcal E_1\ge0$ everywhere on $\Omega$. From Proposition~\ref{prop:kappa-red} with $\kappa=0$, we conclude that $I\ge18\pi$.

If $I=18\pi$, then \eqref{eq:kappa-red} gives $\int_\Omega\mathcal E_1\,d\mu=0$. Continuity gives $\mathcal E_1\equiv0$ on $\Omega$. All three coefficients of $\mathcal E_1$ are positive on $\{0<a<1\}$, and $\mathcal E_1\ge\frac{2-\omega_3}{\omega_3}T$ on $\{a(1-a)=0\}$. Hence $T=0$ and $\Theta=0$ everywhere on $\Omega$.
\end{proof}

In the previous proof, we used the following elementary lemma.

\begin{lemma}\label{lem:d3-profiles}
For $d=3$:
\begin{enumerate}
\item[\textup{(i)}] $w_3$ is concave on $[0,1]$, $|\gamma_3(a)|\le\frac12|1-2a|$ on $[0,1]$, and $0<h_3(a)<\omega_3-1$ on $(0,1/2)$.
\item[\textup{(ii)}] $1-\frac{2}{\omega_3}(1-a)h_3(a)\ge\frac{2-\omega_3}{\omega_3}>0$ and $1+\frac{2}{\omega_3}a\,h_3(a)\ge\frac{2-\omega_3}{\omega_3}>0$ for all $a\in[0,1]$, and
\[
 -\frac{1}{\omega_3}\bigl(h_3'(a)-6\gamma_3(a)h_3(a)\bigr)-6\gamma_3(a)^2>0\qquad\text{for all }a\in(0,1).
\]
\end{enumerate}
\end{lemma}

\begin{proof}
\textup{(i)} By symmetry under $a\leftrightarrow1-a$, let $a\in(0,1/2]$ and set $\xi=\sqrt{a/(1-a)}\in(0,1]$. A direct computation shows that
\[
 w_3''(a)=-\frac{3(\xi^2+1)^5(\xi^2-4\xi+1)^2}{2\xi(1+\xi^3)^4}\le0,\qquad
 \frac{1-2a}{2}-\gamma_3(a)=\frac{\xi(1-\xi)^3}{2(1+\xi^2)(1+\xi^3)}\ge0.
\]
Thus $w_3$ is concave on $[0,1]$ and $|\gamma_3(a)|\le\frac12|1-2a|$, and Lemma~\ref{lem:h-ode} implies $0<h_3<\omega_3-1$ on $(0,1/2)$.

\textup{(ii)} By symmetry, take $a\in[0,1/2]$. Since $0\le h_3\le\omega_3-1$ and $\omega_3<2$, we have $1+\frac{2}{\omega_3}ah_3\ge1$ and $1-\frac{2}{\omega_3}(1-a)h_3\ge\frac{2-\omega_3}{\omega_3}>0$. For the third inequality, at $a=1/2$ \eqref{eq:h-ode} implies $-h_3'(1/2)/\omega_3=6(2-\omega_3)/\omega_3>0$. For $a\in(0,1)$ with $m=\sqrt{a(1-a)}\in(0,1/2)$, we have $|1-2a|=\sqrt{1-4m^2}$ and $a^{3/2}+(1-a)^{3/2}=(1-m)\sqrt{1+2m}$, so
\begin{equation}\label{eq:d3-w-gamma}
 w_3=\frac{1}{(1-m)^2(1+2m)},\qquad |\gamma_3|=\frac{\sqrt{1-4m^2}}{2(1-m)(1+2m)}.
\end{equation}
Restricting to $a\in(0,1/2)$ (so $\gamma_3=|\gamma_3|$) and eliminating $h_3'$ via \eqref{eq:h-ode}, we find
\[
 -\frac{1}{\omega_3}\bigl(h_3'(a)-6\gamma_3(a)h_3(a)\bigr)-6\gamma_3(a)^2=\frac{3w_3}{2\omega_3m^2}\Bigl[(1-m^2)\sqrt{1-4m^2}\,h_3+1-\omega_3(1-2m^2)\Bigr].
\]
In terms of $m$, \eqref{eq:h-ode} reads $\frac{\mathrm d}{\mathrm dm}(m^3h_3)=3m^2(1-4m^2)^{-1/2}(\omega_3-w_3)$, and direct differentiation shows that
\[
 \frac{\mathrm d}{\mathrm dm}\biggl[m^3h_3+\frac{m^3\bigl(1-\omega_3(1-2m^2)\bigr)}{(1-m^2)\sqrt{1-4m^2}}\biggr]
 =\frac{m^4\bigl[\omega_3(1-m^2)(1-4m^2)+6m(1-\omega_3m)\bigr]}{(1-m^2)^2(1-4m^2)^{3/2}}>0
\]
on $(0,1/2)$, since $\omega_3<2$ implies $\omega_3m<1$. Because the bracketed function vanishes at $m=0$, the claimed positivity follows on $(0,1/2)$.
\end{proof}

\subsection{Chiral pairs in dimension four}
In any even dimension $d\ge2$, the chirality operator $\Gamma=i^{d/2}c_1\cdots c_d$ satisfies $\Gamma^*=\Gamma$, $\Gamma^2=I_{N_d}$, and $c_j\Gamma=-\Gamma c_j$, decomposing $\C^{N_d}=\Sigma^+\oplus\Sigma^-$ into orthogonal $\pm1$ eigenspaces of real dimension $N_d$.

\begin{lemma}\label{lem:chiral}
Let $d\ge2$ be even. For $\psi=\psi_++\psi_-\in\X_d\setminus\{0\}$ decomposed along $\C^{N_d}=\Sigma^+\oplus\Sigma^-$, one has $|\psi|^2=|\psi_+|^2+|\psi_-|^2$ and $|D_d\psi|^2=|D_d\psi_+|^2+|D_d\psi_-|^2$ pointwise, and
\begin{equation}\label{eq:chiral-split}
 Q_d(\psi)^2\ge\frac{\norm{D_d\psi_+}_{p_d}^2+\norm{D_d\psi_-}_{p_d}^2}{\norm{\psi_+}_{q_d}^2+\norm{\psi_-}_{q_d}^2}\ge\min_{\psi_\pm\ne0}Q_d(\psi_\pm)^2.
\end{equation}
When both $\psi_+\ne0$ and $\psi_-\ne0$, equality holds throughout if and only if $Q_d(\psi_+)=Q_d(\psi_-)=Q_d(\psi)$ and there exists $\kappa_0>0$ such that $|\psi_+|=\kappa_0|\psi_-|$ and $|D_d\psi_+|=\kappa_0|D_d\psi_-|$ almost everywhere.
\end{lemma}

\begin{proof}
Since $c_j\Gamma=-\Gamma c_j$, the operator $D_d$ exchanges $\Sigma^+$ and $\Sigma^-$, which proves the pointwise norm identities. Because $p_d/2<1<q_d/2$, the reverse Minkowski inequality in $L^{p_d/2}$ and Minkowski's inequality in $L^{q_d/2}$ read
\[
 \norm{D_d\psi}_{p_d}^2\ge\norm{D_d\psi_+}_{p_d}^2+\norm{D_d\psi_-}_{p_d}^2,\qquad
 \norm{\psi}_{q_d}^2\le\norm{\psi_+}_{q_d}^2+\norm{\psi_-}_{q_d}^2,
\]
and \eqref{eq:chiral-split} follows. When both $\psi_\pm\ne0$, equality in the second inequality holds only if $Q_d(\psi_+)=Q_d(\psi_-)$, while strict concavity of $s\mapsto s^{p_d/2}$ and strict convexity of $s\mapsto s^{q_d/2}$ imply $|D_d\psi_+|=\kappa_2|D_d\psi_-|$ and $|\psi_+|=\kappa_1|\psi_-|$ almost everywhere ($\kappa_1,\kappa_2>0$), with $\kappa_1=\kappa_2$ since $Q_d(\psi_+)=Q_d(\psi_-)$. The converse is immediate.
\end{proof}

Now fix $d=4$, so $\dim_\R\Sigma^+=\dim_\R\Sigma^-=4$, $p=8/5$, $q=8/3$, $S_4^4=12\pi^3$, $V_0=128\pi^2/3$, and $\omega_4=9\pi/16$. We call a solution $(U,V)\in\X_4^2$ of \eqref{eq:coupled} \emph{chiral} if $\Gamma U=\pm U$ and $\Gamma V=\mp V$. Since $\omega_4>3/2$, the coefficient $1-\frac{3}{\omega_4}(1-a)h_4(a)$ in $\mathcal E_1$ is negative at $a=0$; we compensate this by choosing $\kappa=1/2$ in Proposition~\ref{prop:kappa-red}.

\begin{proposition}\label{prop:d4-chiral}
Every nontrivial chiral solution $(U,V)\in\X_4^2$ of \eqref{eq:coupled} satisfies $I\ge12\pi^3$. If $I=12\pi^3$, then $T=0$ and $\Theta=0$ everywhere on $\Omega$.
\end{proposition}

\begin{proof}
Wherever $J\ne0$ on $\Omega$, the spinors $c_1J,\dots,c_4J$ lie in the half-spinor space $\Sigma^\pm$ containing $H$ and form a real orthogonal basis of $\Sigma^\pm$ of squared norm $(1-a)\lambda^{-1}$. By Parseval's identity, $|B|^2=a(1-a)$ and $\Theta=|X|^2$ on $\Omega$.

We claim that $\mathcal E_0\ge0$ and $\mathcal E_1+\frac12\mathcal E_0\ge0$ pointwise on $\Omega$, with
\begin{equation}\label{eq:d4-comb-pos}
 4a(1-a)\Bigl(\mathcal E_1+\frac12\mathcal E_0\Bigr)\ge\frac{2-\omega_4}{\omega_4}\,|X|^2\qquad\text{on }\{0<a<1\}.
\end{equation}
On $\{a(1-a)=0\}$, we have $X=0$, so $\mathcal E_0=T\ge0$, and $|h_4|\le\omega_4-1$ (Lemma~\ref{lem:d4-profiles}\textup{(i)} below) implies $\frac23(\mathcal E_1+\frac12\mathcal E_0)\ge\frac{2-\omega_4}{\omega_4}T\ge0$. On $\{0<a<1\}$, by $(U,V)\leftrightarrow(V,U)$ symmetry we may assume $a\in(0,1/2]$, where $0\le\gamma_4\le(1-2a)/2$ by Lemma~\ref{lem:d4-profiles}\textup{(i)}. Thus
\[
 2\bigl(1-32a(1-a)\gamma_4^2\bigr)\ge2-16a(1-a)\bigl(1-4a(1-a)\bigr)=1+\bigl(1-8a(1-a)\bigr)^2\ge1,
\]
and Lemma~\ref{lem:twistor-lower} at $\tau=0$ implies $6a(1-a)\mathcal E_0\ge|X|^2\ge0$.

Next, set $\tau=(2/\omega_4)h_4(a)$ and $\tau_0=(2/\omega_4)(\omega_4-1)$, so $0\le\tau\le\tau_0<1$ with $1-\tau_0=(2-\omega_4)/\omega_4>0$. Applying \eqref{eq:twistor-weighted} at $d=4$ to $\frac23\bigl(\mathcal E_1+\frac12\mathcal E_0\bigr)$, we find
\begin{align*}
 4a(1-a)\Bigl(\mathcal E_1+\frac12\mathcal E_0\Bigr)
 &\ge2\Bigl(1-\bigl(1-2a-8a(1-a)\gamma_4\bigr)\tau-32a(1-a)\gamma_4^2\Bigr)|X|^2\\
 &\quad+6a(1-a)\Bigl(-\frac{h_4'}{\omega_4}\Bigr)|X|^2.
\end{align*}
Since $\gamma_4\le(1-2a)/2$, we have $1-2a-8a(1-a)\gamma_4\ge(1-2a)^3\ge0$, and we may replace $\tau$ by $\tau_0$. Moreover, by Lemma~\ref{lem:d4-profiles}\textup{(ii)} below and $\omega_4<2$,
\[
 6a(1-a)\Bigl(-\frac{h_4'}{\omega_4}\Bigr)=\frac{3a(1-a)\tau_0}{\omega_4-1}(-h_4')\ge4a\tau_0\Bigl[w_4(a)\bigl(1+4(1-2a)\gamma_4\bigr)-2\Bigr].
\]
Using $w_4(a)\ge1+2a$ (by concavity on $[0,1/2]$) and $2\gamma_4^2\le(1-2a)\gamma_4$, we bound the coefficient of $\tau_0|X|^2$ from below by
\[
 4a\Bigl[w_4(a)\bigl(1+4(1-2a)\gamma_4\bigr)-1+4(1-a)\gamma_4-16(1-a)\gamma_4^2\Bigr]\ge8a^2\Bigl[1+2\bigl(1+4(1-2a)\bigr)\gamma_4\Bigr]>0.
\]
Together with $2(1-\tau_0)(1-32a(1-a)\gamma_4^2)\ge1-\tau_0=(2-\omega_4)/\omega_4$, this proves \eqref{eq:d4-comb-pos}.

Assume now that $I\le12\pi^3$. By Proposition~\ref{prop:kappa-red} at $d=4$ and $\kappa=1/2$, we have $1\le s^4\le\omega_4<2<4$ (so $1\le s^2<2$) and
\[
 \frac34\Bigl(\frac{I}{12\pi^3}-1\Bigr)V_0\ge\frac38(s^2-1)(2-s^2)V_0+\int_\Omega\Bigl(\mathcal E_1+\frac12\mathcal E_0\Bigr)\,d\mu\ge0.
\]
Thus $I=12\pi^3$ and $\int_\Omega\bigl(\mathcal E_1+\frac12\mathcal E_0\bigr)\,d\mu=0$. By continuity, $\mathcal E_1+\frac12\mathcal E_0\equiv0$ on $\Omega$. In view of \eqref{eq:d4-comb-pos}, we have $X=0$ on $\{0<a<1\}$ and hence on $\Omega$; since $1-\tau_0>0$, it follows that $T=0$ everywhere on $\Omega$.
\end{proof}

In the previous proof, we used the following bounds for $d=4$.

\begin{lemma}\label{lem:d4-profiles}
For $d=4$:
\begin{enumerate}
\item[\textup{(i)}] $w_4$ is strictly concave on $[0,1]$, $0\le\gamma_4(a)\le\frac12(1-2a)$ on $[0,1/2]$, and $0<h_4(a)<\omega_4-1$ on $(0,1/2)$.
\item[\textup{(ii)}] For all $a\in(0,1/2]$,
\begin{equation}\label{eq:h4-prime}
 -h_4'(a)\ge\frac{4}{3(1-a)^2}\Bigl[w_4(a)\bigl(1+4(1-2a)\gamma_4(a)\bigr)-\omega_4\Bigr]\ge\frac{4}{3(1-a)^2}(2-\omega_4)>0.
\end{equation}
\end{enumerate}
\end{lemma}

\begin{proof}
\textup{(i)} For $a\in(0,1/2]$, put $\xi=(a/(1-a))^{1/3}\in(0,1]$. Direct differentiation shows that
\[
 w_4''(a)=-\frac{4(1+\xi^3)^6\bigl[(\xi^3-2\xi^2-2\xi+1)^2+\xi(\xi-1)^2\bigl(4(\xi-1)^2+\xi\bigr)\bigr]}{3\xi^2(1+\xi^4)^5}<0
\]
and
\[
 \frac{1-2a}{2}-\gamma_4(a)=\frac{\xi(1-\xi)\bigl[(1-\xi^2)^2+\xi(1-\xi)^2+2\xi^2\bigr]}{2(1+\xi^3)(1+\xi^4)}\ge0.
\]
By symmetry, $w_4$ is strictly concave on $[0,1]$, and $0<h_4<\omega_4-1$ on $(0,1/2)$ by Lemma~\ref{lem:h-ode}.

\textup{(ii)} On $[0,a]$ ($0<a\le1/2$), we have $1-s\ge1-a$, $w_4(a)-w_4(s)\ge8\gamma_4(a)w_4(a)(a-s)$, and $\int_0^a s(a-s)\,ds=a^3/6$. Inserting these bounds into \eqref{eq:hd-prime-int} at $d=4$ and recalling from the proof of Lemma~\ref{lem:h-ode} that $w_4(a)[1+4(1-2a)\gamma_4(a)]\ge2>\omega_4$, we obtain \eqref{eq:h4-prime}.
\end{proof}

\begin{remark}
Section~\ref{sec:euclidean} and Proposition~\ref{prop:kappa-red} hold for every $d\ge3$. When $d\ge5$, however, $\dim_\R\C^{N_d}>d+1$ for odd $d$ ($\dim_\R\C^{N_5}=8>6$) and $\dim_\R\Sigma^\pm>d$ for even $d$ ($\dim_\R\Sigma^\pm=8>6$ at $d=6$), so the dimension count no longer forces $a(1-a)-|B|^2$ to equal $\lambda^2\ip{H}{J}^2$ in odd dimensions or to vanish for chiral pairs in even dimensions.
\end{remark}

\section{Cases of equality}

In this section we classify all nontrivial solutions of \eqref{eq:coupled} on $\R^d$ ($d\ge3$) satisfying $T=0$ and $\Theta=0$ on $\Omega$, and complete the proofs of Theorems~\ref{thm} and~\ref{thm:pairs} and Corollary~\ref{cor}.

For $x\in\R^d\setminus\{0\}$, write $y=x/|x|^2$. Since $|W_{\mathcal K}(x)|=|x|^{1-d}|W(y)|$, the Kelvin transform \eqref{eq:kelvin} acts on $\lambda$, $a$, and $(H,J)$ off the origin by
\begin{equation}\label{eq:kelvin-HJ}
 \begin{gathered}
  \lambda(x)\longmapsto|x|^{-2}\lambda(y),\qquad a(x)\longmapsto a(y),\\
  (H,J)(x)\longmapsto\bigl(c(x)H(y),\,-c(x)J(y)\bigr).
 \end{gathered}
\end{equation}
Setting $R_{jk}(x)=\delta_{jk}-2x_jx_k/|x|^2$, we see from $c(x)c_jc(x)=|x|^2\sum_{k=1}^dR_{jk}(x)c_k$ and $\partial_jy_k=|x|^{-2}R_{jk}(x)$ that \eqref{eq:kelvin} maps $B(y)$ and $X(y)$ to $R(x)B(y)$ and $R(x)X(y)$, respectively, and gives $T(y)$ and $\Theta(y)$ as the transformed scalar quantities at $x$. Translations and dilations $(H,J)(x)\mapsto\sigma^{1/2}(H,J)(x/\sigma)$ ($\sigma>0$) also preserve $T$ and $\Theta$ under pullback.

\begin{proposition}\label{prop:killing}
Let $d\ge3$. A nontrivial solution $(U,V)\in\X_d^2$ of \eqref{eq:coupled} satisfies $T=0$ and $\Theta=0$ everywhere on $\Omega$ if and only if $(U,V)\in\Or_d$. Moreover, every pair $(U,V)\in\Or_d$ satisfies $\Omega=\R^d$ and $\ip{U}{V}=0$ on $\R^d$, and $\{\kappa U:\kappa>0,\ (U,V)\in\Or_d\}$ coincides with the family \eqref{eq:optimizers}.
\end{proposition}

\begin{proof}
\emph{Step~1: Twistor spinors and $\Omega=\R^d$.} Assume $T=0$ and $\Theta=0$ on $\Omega$, so that $\T H=\T J=0$, $X=0$, and $|B|^2=a(1-a)$ on $\Omega$. Let $O$ be a connected component of $\Omega$.

We first recall (cf.\ \cite[Theorem~10]{FL}) why any $\varphi\in C^2(O;\C^{N_d})$ solving $\partial_j\varphi=-\frac1dc_jD_d\varphi$ ($1\le j\le d$) has the form $\varphi(x)=A+c(x)C$ on $O$. Differentiating the $j$th equation in $x_j$ and summing over $j$, we find $\Delta\varphi=\frac1d\Delta\varphi$, so $\Delta\varphi=0$ and $\varphi\in C^\infty(O;\C^{N_d})$. Differentiating the twistor equation in $x_k$, multiplying on the left by $c_j$, and summing over $j\ne k$ using $D_d^2\varphi=0$, we obtain
\[
 c_j\partial_kD_d\varphi=c_k\partial_jD_d\varphi,\qquad (d-2)\partial_kD_d\varphi=0\qquad(1\le j,k\le d).
\]
Since $d\ge3$, $D_d\varphi$ is constant on $O$, say $D_d\varphi=-dC$, whence $\nabla(\varphi-c(x)C)=0$ and $\varphi(x)=A+c(x)C$ on $O$.

Consequently, $H(x)=A_1+c(x)C_1$ and $J(x)=A_2+c(x)C_2$ on $O$. Along any sequence in $O$ approaching a finite boundary point $x_*\in\partial O$, the polynomial $\lambda(x)^{-1}=|H(x)|^2+|J(x)|^2$ remains bounded, whereas $\lambda(x)\le\rho(x)^\nu\to0$. Thus $\partial O=\emptyset$ and $\Omega=\R^d$.

By \eqref{eq:bessel-B} and $|B|^2=a(1-a)$, we have $\ip{H}{J}=0$ on $\R^d$. The Clifford relations and the definition of $B$ give
\[
 \lambda|aJ+c(B)H|^2=a\bigl(a(1-a)-|B|^2\bigr),\qquad
 \lambda|(1-a)H-c(B)J|^2=(1-a)\bigl(a(1-a)-|B|^2\bigr).
\]
Both right-hand sides vanish, so
\[
 aJ=-c(B)H,\qquad (1-a)H=c(B)J\qquad\text{on }\R^d.
\]
Inserting $\nabla a=-\frac2d\lambda B$ and these two identities into \eqref{eq:D-Phi} for $\Phi_1=\lambda^{d/2}H$ and $\Phi_2=\lambda^{d/2}J$, we see that $\gamma(a)$ cancels and
\[
 D_dH+\frac{d}{2}\lambda^{-1}c(\nabla\lambda)H=\lambda J,\qquad
 D_dJ+\frac{d}{2}\lambda^{-1}c(\nabla\lambda)J=\lambda H.
\]
Put $L_\pm=H\pm J=A_\pm+c(x)C_\pm$ on $\R^d$. Because $\ip{H}{J}=0$, we have
\[
 \lambda(x)^{-1}=|H(x)|^2+|J(x)|^2=|L_\pm(x)|^2=\alpha_0+2v_0\cdot x+\alpha_2|x|^2
\]
with $\alpha_0=|A_\pm|^2=\lambda(0)^{-1}>0$, $\alpha_2=|C_\pm|^2\ge0$, and $v_0\in\R^d$. Using $D_dL_\pm=-dC_\pm$ and $-\lambda^{-2}\nabla\lambda=2(v_0+\alpha_2x)$, we add and subtract the equations for $H$ and $J$ to obtain
\begin{equation}\label{eq:Lpm-alg}
 \lambda(x)^{-1}C_\pm=\Bigl(\mp\frac1d-c(v_0)-\alpha_2c(x)\Bigr)\bigl(A_\pm+c(x)C_\pm\bigr).
\end{equation}

\emph{Step~2: Reduction to $(U_\eta,V_\eta)$.} At $x=0$, \eqref{eq:Lpm-alg} reads $\alpha_0C_+=-(1/d+c(v_0))A_+$; taking squared norms with $\ip{A_+}{c(v_0)A_+}=0$ shows that $\alpha_0\alpha_2=|v_0|^2+d^{-2}>0$. Since $H$ is affine and $\lambda(x)^{-1}$ is quadratic with $\alpha_2>0$, the ratio $a(x)=\lambda(x)|H(x)|^2$ has a finite limit $a_\infty\in[0,1]$ as $|x|\to\infty$, while $X=0$ and $|B|^2=a(1-a)$ imply $|\nabla a|^2=\frac{4}{d^2}\lambda^2a(1-a)$ on $\R^d$.

We use the conformal symmetries to put a zero of $a$ at infinity. If $a$ attains its minimum at $x_{\min}\in\R^d$, the gradient identity gives $a(x_{\min})\in\{0,1\}$. Since $\int_{\R^d}|V|^q\,dx>0$, the minimum is zero. Translate $x_{\min}$ to $0$ and apply \eqref{eq:kelvin-HJ} to send this zero to infinity. If the minimum is attained only at infinity, apply the same finite-minimum argument to the Kelvin image, whose ratio extends across $0$ by Step~1 and has its minimum there. This proves $a_\infty=0$ for the original pair, which we retain in this case.

After this transformation, $\lim_{|x|\to\infty}a(x)=0$. Completing the square in $\lambda(x)^{-1}$ and applying a translation and critical dilation, we may assume that $\lambda(x)=d/(1+|x|^2)$ while keeping $\lim_{|x|\to\infty}a(x)=0$. Then $\alpha_0=\alpha_2=1/d$ and $v_0=0$, so \eqref{eq:Lpm-alg} at $x=0$ implies $C_\pm=\mp A_\pm$ with $|A_\pm|^2=1/d$, and
\[
 a(x)=\frac{|L_+(x)+L_-(x)|^2}{4\lambda(x)^{-1}}=\frac{d\,\bigl|(A_++A_-)-c(x)(A_+-A_-)\bigr|^2}{4(1+|x|^2)}\longrightarrow\frac{d}{4}|A_+-A_-|^2
\]
as $|x|\to\infty$. Hence $A_+=A_-$. With $\eta=\sqrt d\,A_+$ ($|\eta|=1$), we find $(H,J)(x)=d^{-1/2}(\eta,-c(x)\eta)$, and $(U,V)=(|H|^q+|J|^q)^{-(d-1)/2}(H,J)=(U_\eta,V_\eta)$, so the original pair lies in $\Or_d$. Conversely, $(U_\eta,V_\eta)$ satisfies $T=0$, $\Theta=0$, $\Omega=\R^d$, and $\ip{U_\eta}{V_\eta}=0$, and conformal invariance extends these properties to $\Or_d$.

\emph{Step~3: Identification with \eqref{eq:optimizers}.} For $x\in\R^d\setminus\{0\}$, $s\in\R$, and $\xi\in\R^d$, let $\widehat\xi=\xi/|\xi|$ if $\xi\ne0$ and let $\widehat\xi$ be any unit vector if $\xi=0$; then
\begin{equation}\label{eq:inv-norm}
 |x|\,\Bigl|s\frac{x}{|x|^2}+\xi\Bigr|=\bigl|\,|\xi|x+s\widehat\xi\,\bigr|.
\end{equation}
Translations and dilations of $f_d(|x|)\eta$ ($\eta\ne0$) form the first family in \eqref{eq:optimizers}. Translating any member of the second family by $-x_1$, we obtain $\psi_0(y)=(|y|^q+|\alpha y+\xi_0|^q)^{-(d-1)/2}c(y)\eta$ with $\xi_0=\alpha x_1+\beta\ne0$, which by \eqref{eq:inv-norm} coincides with $(\psi_1)_{\mathcal K}$ for
\[
 \psi_1(x)=\frac{\eta}{\bigl(1+|\xi_0|^q|x+\alpha|\xi_0|^{-2}\xi_0|^q\bigr)^{(d-1)/2}},
\]
so every spinor in \eqref{eq:optimizers} belongs to $\{\kappa U:\kappa>0,\ (U,V)\in\Or_d\}$. Conversely, \eqref{eq:optimizers} contains $f_d(|x|)\eta$ and is invariant under translations, dilations, and $\psi\mapsto\psi_{\mathcal K}$: by \eqref{eq:inv-norm}, $\psi\mapsto\psi_{\mathcal K}$ exchanges the first family with the $x_1=0$ case of the second, whereas for $x_1\ne0$ in the second family, with $\widetilde x_1=x_1/|x_1|^2$ we have $c(x)c(x/|x|^2-x_1)\eta=-c(x-\widetilde x_1)c(x_1)\eta$, so $\psi_{\mathcal K}$ belongs to the second family with parameters $(\widetilde x_1,\widetilde\alpha,\widetilde\beta)=\bigl(\widetilde x_1,\,|\beta|/|x_1|,\,(\alpha/|x_1|)\widehat\beta\bigr)$ satisfying $|\widetilde\alpha\widetilde x_1+\widetilde\beta|=|\alpha x_1+\beta|/|x_1|^2>0$.
\end{proof}

It remains to remove the chirality assumption for $d=4$.

\begin{lemma}\label{lem:d4-nonchiral}
Every nontrivial solution $(U,V)\in\X_4^2$ of \eqref{eq:coupled} satisfies $I\ge12\pi^3$, with equality if and only if $(U,V)\in\Or_4$.
\end{lemma}

\begin{proof}
By Proposition~\ref{prop:attain} and Lemma~\ref{lem:chiral}, $Q_4$ has a minimizer on $\X_4\setminus\{0\}$ taking values in one half-spinor space $\Sigma^\pm$. Since $z'\mapsto|z'|^{p_4-2}z'$ preserves $\Sigma^\pm$ and $D_4$ exchanges $\Sigma^+$ and $\Sigma^-$, its normalization \eqref{eq:normalize} is a chiral solution of \eqref{eq:coupled}, so Proposition~\ref{prop:d4-chiral} and Lemma~\ref{lem:trial} imply $\inf_{\X_4\setminus\{0\}}Q_4=S_4=(12\pi^3)^{1/4}$. By \eqref{eq:equal-norm}, every nontrivial solution $(U,V)\in\X_4^2$ of \eqref{eq:coupled} therefore satisfies $I=Q_4(U)^4\ge12\pi^3$.

Suppose now that $I=12\pi^3$, so $Q_4(U)=S_4$ and $\norm{U}_{8/3}^{8/3}=12\pi^3$, and write $U=\psi_++\psi_-$ with $\psi_\pm$ taking values in $\Sigma^\pm$. If $\psi_+=0$ or $\psi_-=0$, then $(U,V)$ is chiral and belongs to $\Or_4$ by Propositions~\ref{prop:d4-chiral} and~\ref{prop:killing}.

Suppose both $\psi_\pm\ne0$. The equality conditions in Lemma~\ref{lem:chiral} give $Q_4(\psi_+)=Q_4(\psi_-)=S_4$ and $\theta\in(0,\pi/2)$ such that $|\psi_+|=(\cos\theta)|U|$, $|\psi_-|=(\sin\theta)|U|$, and $|D_4\psi_+|=(\cot\theta)|D_4\psi_-|$ almost everywhere. Set $U_+=(\cos\theta)^{-1}\psi_+$, $U_-=(\sin\theta)^{-1}\psi_-$, and $V_\pm=|D_4U_\pm|^{p_4-2}D_4U_\pm$. Each $(U_\pm,V_\pm)$ is a chiral solution of \eqref{eq:coupled} with $\int_{\R^4}|U_\pm|^{8/3}\,dx=12\pi^3$. By Propositions~\ref{prop:d4-chiral} and~\ref{prop:killing}, $(U_\pm,V_\pm)\in\Or_4$, with $|U_+|=|U_-|=|U|$, $|V_+|=|V_-|$, and $\rho>0$ on $\R^4$.

The pointwise proportionality above gives $|D_4U_+|=|D_4U_-|$. Orthogonality of the chiral components then gives $|D_4U|=|D_4U_\pm|=|V_\pm|^{q_4-1}$. Hence
\[
 U=\cos\theta\,U_++\sin\theta\,U_-,\qquad V=|D_4U|^{p_4-2}D_4U=\cos\theta\,V_++\sin\theta\,V_-.
\]
Thus $(U,V)$ has the same $\rho>0$, $a$, and $\lambda$ on $\R^4$ as $(U_\pm,V_\pm)$. By linearity of $\T$, we have $\T H=\T J=0$ ($T=0$). Since $\max\{1+2a\gamma_4,\,1-2(1-a)\gamma_4\}\ge1$, \eqref{eq:TH-H} gives $X=0$. It follows that $B=-2\lambda^{-1}\nabla a=B_+$ and $|B|^2=|B_+|^2=a(1-a)$ ($\Theta=0$). Proposition~\ref{prop:killing} now gives $(U,V)\in\Or_4$. The converse follows from Lemmas~\ref{lem:trial} and~\ref{lem:conf}.
\end{proof}

\begin{proof}[Proof of Theorem~\ref{thm:pairs}]
For $d=3$, this follows from Propositions~\ref{prop:d3} and~\ref{prop:killing} together with Lemmas~\ref{lem:trial} and~\ref{lem:conf}. For $d=4$, this is Lemma~\ref{lem:d4-nonchiral}.
\end{proof}

\begin{proof}[Proof of Theorem~\ref{thm}]
Section~\ref{sec:d2} proves Theorem~\ref{thm} for $d=2$. Let $d\in\{3,4\}$. By Proposition~\ref{prop:attain}, $Q_d$ has a global minimizer $\psi\in\X_d\setminus\{0\}$; its normalization \eqref{eq:normalize} satisfies $I=Q_d(\psi)^d\ge S_d^d$ by Theorem~\ref{thm:pairs}, so $\inf Q_d=S_d$ by Lemma~\ref{lem:trial}. If $\psi\in\X_d\setminus\{0\}$ attains $Q_d(\psi)=S_d$, then its normalization $(U,V)$ in \eqref{eq:normalize} has $I=S_d^d$, so $(U,V)\in\Or_d$ by Theorem~\ref{thm:pairs} and $\psi$ belongs to $\{\kappa U:\kappa>0,\ (U,V)\in\Or_d\}$. Conversely, every element of this set attains $Q_d(\psi)=S_d$ by Lemmas~\ref{lem:crit-norm}, \ref{lem:trial}, and~\ref{lem:conf}, and Proposition~\ref{prop:killing} identifies it with \eqref{eq:optimizers}.
\end{proof}

\begin{proof}[Proof of Corollary~\ref{cor}]
Fix a unit spinor $\eta_0\in\C^{N_d}$. Given $u\in L^1_{\mathrm{loc}}(\R^d;\R)$ such that $\nabla u\in L^{q_d}(\R^d;\R^d)$ and $\Delta u\in L^{p_d}(\R^d)$, mollification shows that $\psi=c(\nabla u)\eta_0\in\X_d$ satisfies $D_d\psi=-(\Delta u)\eta_0$ distributionally with $|\psi|=|\nabla u|$ and $|D_d\psi|=|\Delta u|$ pointwise. Inequality \eqref{eq:riesz-grad} therefore follows from Theorem~\ref{thm}.

Suppose $\nabla u\not\equiv0$ attains equality in \eqref{eq:riesz-grad}, so $Q_d(\psi)=S_d$. Set
\[
 \Phi=|D_d\psi|^{p_d-2}D_d\psi=-|\Delta u|^{p_d-2}(\Delta u)\,\eta_0\in L^{q_d}(\R^d;\C^{N_d}),
\]
equal to $0$ where $\Delta u=0$. By \eqref{eq:EL}, $D_d\Phi=\Lambda|\psi|^{q_d-2}\psi$ distributionally on $\R^d$, with $\norm{\Phi}_{q_d}=\norm{D_d\psi}_{p_d}^{p_d-1}$ and $\norm{D_d\Phi}_{p_d}=\Lambda\norm{\psi}_{q_d}^{q_d-1}$, so $\Phi\in\X_d\setminus\{0\}$ also attains $Q_d(\Phi)=S_d$.

By Theorem~\ref{thm}, $\Phi$ agrees almost everywhere with a spinor in \eqref{eq:optimizers}. Moreover, $\Phi$ takes values in $\R\eta_0$ almost everywhere. Since $y\mapsto c(y)\eta$ is injective on $\R^d$ and $d\ge2$, no member of the second family in \eqref{eq:optimizers} can take values in $\R\eta_0$ almost everywhere. Hence $\Phi(x)=\kappa\eta_0(\sigma^{q_d}+|x-x_0|^{q_d})^{-(d-1)/2}$ almost everywhere for some $\kappa\in\R\setminus\{0\}$, $\sigma>0$, and $x_0\in\R^d$. Differentiation gives
\[
 \Lambda|\psi(x)|^{q_d-2}\psi(x)=D_d\Phi(x)=-\frac{d\kappa\,|x-x_0|^{q_d-2}c(x-x_0)\eta_0}{\bigl(\sigma^{q_d}+|x-x_0|^{q_d}\bigr)^{(d+1)/2}}\qquad\text{a.e.\ on }\R^d.
\]
Applying $z'\mapsto|z'|^{p_d-2}z'$ and comparing with $\psi=c(\nabla u)\eta_0$, we obtain \eqref{eq:riesz-eq} with $C=-(d/\Lambda)^{p_d-1}|\kappa|^{p_d-2}\kappa$. Conversely, for any $u_0\in\R$, $C\in\R\setminus\{0\}$, $\sigma>0$, and $x_0\in\R^d$, the radial primitive
\[
 u(x)=u_0+\int_0^{|x-x_0|}\frac{Cr}{\bigl(\sigma^{q_d}+r^{q_d}\bigr)^{(d-1)/2}}\,dr\in C^4(\R^d)
\]
satisfies \eqref{eq:riesz-eq} and $\Delta u(x)=Cd\sigma^{q_d}(\sigma^{q_d}+|x-x_0|^{q_d})^{-(d+1)/2}\in L^{p_d}(\R^d)$, and attains equality in \eqref{eq:riesz-grad} since $\psi=c(\nabla u)\eta_0$ belongs to \eqref{eq:optimizers}.
\end{proof}

\end{document}